\documentstyle[12pt]{article}
\input mssymb
\newtheorem{theorem}{Theorem}[section]
\newtheorem{claim}[theorem]{Claim}

\newtheorem{defn}[theorem]{Definition}
\newtheorem{eple}[theorem]{Example}

\newtheorem{lemma}[theorem]{Lemma}
\newtheorem{prop}[theorem]{Proposition}

\newenvironment{example}{\begin{eple}\em}{\end{eple}}
\newsavebox{\indbin}
\savebox{\indbin}{\begin{picture}(0,0)
\newlength{\gnu}
\settowidth{\gnu}{$\smile$}
\setlength{\unitlength}{.5\gnu}
\put(-1,-.65){$\smile$}
\put(-.25,.1){$|$}
\end{picture}}

\newcommand{\be}{\begin{enumerate}}
\newcommand{\bd}{\begin{defn}}
\newcommand{\bt}{\begin{theorem}}
\newcommand{\bl}{\begin{lemma}}
\newcommand{\ee}{\end{enumerate}}
\newcommand{\ed}{\end{defn}}
\newcommand{\et}{\end{theorem}}
\newcommand{\el}{\end{lemma}}

\newcommand{\la}{\langle}
\newcommand{\ra}{\rangle}
\newcommand{\sub}{\subseteq}

\newcommand{\sm}{\setminus}

\newcommand{\proof}{{\bf Proof:}\quad}

\newcommand{\Cl}{{\cal L}}

\newcommand{\th}{{\mbox{\scriptsize th}}}
\newcommand{\gen}{\mbox{\rm Gen}}
\newcommand{\col}{{col}}

\newcommand{\otp}{\mbox{otp}\,}

\newcommand{\res}{\!\!\restriction\!}
\newcommand{\CA}{{\cal A}}

\newcommand{\CD}{{\cal D}}

\newcommand{\CG}{{\cal G}}

\newcommand{\CL}{{\cal L}}

\newcommand{\CP}{{\cal P}}
\newcommand{\CR}{{\cal R}}

\def\proves{\vdash}

\def\ad{{\alpha\delta}}

\def\bigspace{\vskip.1in }
\def\cof{\mbox{cf}\,}
\def\containing{\supseteq} \def\contains{\supseteq}

\def\dom{{\rm dom}\,}
\def\fakebf#1{\rlap{$#1$}\kern.2pt\rlap{$#1$}\kern.2pt #1}

\def\hgt{{\rm ht}\,}
\def\implies{\rightarrow}
\def\includedin{\subseteq}
\def\intersect{\cap}
\def\iso{\simeq}
\def\ld{{<\!\delta}}
\def\len{{\rm len}\,}
\def\lesd{\le_{\rm sd}}
\def\ll{{<\!\lambda}}

\newtheorem{notation}[theorem]{Notation}
\def\nl{\mbox{}\newline}
\def\page{\vfill \break}
\def\pp{{\bf p}}
\def\qq{{\bf q}}
\def\restr{\!\!\restriction\!}

\def\setoff#1\par{\medskip\noindent{\bf #1}\par\medskip}
\def\union{\cup}
\def\up#1{{}^{#1}}
\def\with{\leftrightarrow}
\def\CPll{\CP_{<\lambda}}
\def\Diag{\Delta}
\def\Dl{{\mbox{\rm Dl}}}
\def\Gh{\hat G}
\def\Gbar{\bar  G}
\def\Intersection{\bigcap}
\def\Pp{{\Bbb P}}
\def\Union{\bigcup}
\def\qed{\hfill\mbox{$\Box$}}
\def\lam{\lambda}

\title{Models with second order properties V:\\ A general principle}
\author{S. Shelah \and C. Laflamme \and B. Hart\thanks{This paper is
based on lectures given by the first author while visiting the
University of Michigan.  The research was partially supported by the
BSF, the NSF and the NSERC.}} 
\begin{document} \maketitle
\bibliographystyle{plain}

\noindent {\bf Abstract}.  \quad  We present 
a general framework for carrying out the
constructions in [2-10] and others of
the same type. The unifying factor is a
combinatorial principle which we present
in terms of a game in which the first
player challenges the second player to
carry out constructions which would be
much easier in a generic extension of
the universe, and the second player
cheats with the aid of $\Diamond$.
\S1
contains an axiomatic framework suitable
for the description of a number of
related constructions, and the statement
of the main theorem \ref{main} in terms
of this framework. In \S2 we illustrate
the use of our combinatorial principle.
The proof of the main result is then
carried out in \S\S3-5.
\page

{\noindent\bf Contents\bigskip}

\noindent\S1. 
Uniform partial orders.\medskip

We describe a class of partial orderings
associated with attempts to manufacture
an object of size $\lambda^+$ from
approximations of size less than
$\lambda$.  We also introduce some
related notions motivated by the forcing
method.  The underlying idea is that a
sufficiently generic filter on the given
partial ordering should give rise to the
desired object of size $\lambda^+$.

We describe a game for two players, in
which the first player imposes
genericity requirements on a
construction, and the second player
constructs an object which meets the
specified requirements. The main theorem
(\ref{main}) is that under certain
combinatorial conditions the second
player has a winning strategy for this
game.
\medskip

\noindent\S2. 
Illustrative application.\medskip

We illustrate the content of our general
principle with an example.  We show
the completeness of the logic $\CL^{<\omega}$, defined by Magidor and Malitz
\cite{MagMal} 
for the $\lambda^+$-interpretation
assuming the combinatorial principles $\Dl_\lambda$ and
$\Diamond_{\lambda^+}$. 

\medskip

\noindent\S3. 
Commitments.\medskip

We give a preliminary sketch of the
proof of Theorem
\ref{main}. 
 We then introduce the notion of 
``basic data'' which is a collection of combinatorial
objects derived from $\Dl_\lambda$ and an object called a commitment
describing the main features of
the second player's strategy in a given
play of the genericity game.  We state
the main results concerning commitments,
and show how Theorem \ref{main} follows
from these results.
\medskip

\noindent\S4. 
Proofs.
\medskip

We prove the propositions stated in \S3 except we defer the proof of
Propositions \ref{basestep} and
\ref{successor} to section
\S5. 
 We use $\Dl_\lambda$ to show that a suitable collection of
``basic data'' exists.  Then we
verify some continuity properties
applying to our strategy at limit
ordinals.
\medskip

\noindent 
\S5. Proof of Proposition \ref{successor}.
\medskip

We prove
Proposition \ref{successor} as well as Proposition \ref{basestep}.
\bigskip

{\noindent\bf Notation}

 If $(A_\alpha : \alpha < \delta)$ is an
increasing sequence of sets we write
$A_{<\delta}$ for $\Union_{\alpha <
\delta} 
A_\alpha$.  (Note the exception arising
in lemma
\ref{mingen}.)

Throughout the paper, $\lambda$ is a
cardinal such that $\lambda^{<\lambda} =
\lambda$.

$\CPll(A) = \{B \sub A: |B| <
\lambda\}$.
     
$\otp(u)$ will mean the order type of
$u$.  Trees are well-founded, and if $T$
is a tree, $\eta\in T$, we write
$\len(\eta)$ for $\otp\{\nu\in T:\nu <
\eta\}$ (the level at which $\eta$
occurs in $T$).

\noindent {\bf Acknowledgements}

We would like to express our deepest gratitude to Greg Cherlin without
whose insistence this paper might not have been completed this century.

\vfill 
\break
    
\section{Uniform 
partial orders}

We will present an axiomatic
framework for the construction of
objects of size $\lambda^+$ from
approximations of size less than
$\lambda$, under suitable set
theoretical hypotheses.  The basic idea
is that we are constructing objects
which can fairly easily be forced to
exist in a generic extension, and we
replace the forcing construction by the
explicit construction of a sufficiently
generic object in the ground model.

We begin with the description of the
class of partial orderings to which our
methods apply.  Our idea is that an
``approximation'' to the desired final
object is built from a set of ordinals
$u\includedin
\lambda^+$ 
of size less than $\lambda$.
Furthermore, though there will be many
such sets $u$, there will be at most
$\lambda$ constructions applicable to an
arbitrary set $u$.
We do not axiomatize the notion of a
``construction'' in any detail; we
merely assume that the approximations
can be coded by pairs $(\alpha, u)$,
where $\alpha<\lambda$ is to be thought
of as a code for the particular
construction applied to $u$.  An
additional feature, suggested by the
intuition just described, is captured in
the ``indiscernibility'' condition
below, which is a critical feature of
the situation -- though trivially true
in any foreseeable application.

\begin{defn}

A {\em standard $\lambda^+$-uniform}
partial order is a partial order $\le$
defined on a subset $\Pp$ of $\lambda
\times 
\CP_{<
\lambda}(\lambda^+)$ 
satisfying the following conditions,
where for $p=(\alpha,u)$ in $\Pp$ we
write $\dom p=u$, and call $u$ the {\em
domain} of $p$.
\be \item If $p\le q$ then 
$\dom p\includedin \dom q$.
\item 
For all $p,q,r \in \Pp$ with $p,q \le r$
there is $r' \in \Pp$ so that $p,q \le
r' \le r$ and $\dom r'=\dom p \union
\dom q$.
\item 
If $(p_i)_{i < \delta} $ is an
increasing sequence of length less than
$\lambda$, then it has a least upper
bound $q$, with domain \
$\bigcup_{i<\delta} \dom p_i$; we will
write $q=\bigcup_{i<\delta} p_i$, or
more succinctly: $q=p_{<\delta}$.
\item 
For all $p \in \Pp$ and $\alpha <
\lambda^+$ there exists a $q \in \Pp$
with $q \le p$ and $\dom q = \dom p
\intersect \alpha$; furthermore, there
is a unique maximal such $q$, for which
we write $q=p \restr \alpha$.
\item 
For limit ordinals $\delta$, $p \restr
\delta =
\bigcup_{\alpha 
< \delta} p \restr \alpha$.
\item 
If $(p_i)_{i < \delta}$ is an increasing
sequence of length less than $\lambda$,
then \[(\Union_{i < \delta}
p_i)\res\alpha\ =
\ 
\Union_{i < \delta} 
(p_i \restr  \alpha).\]
\item 
{\em (Indiscernibility)} If $p=(\alpha,
v) \in \Pp$ and $h:v \rightarrow
v'\includedin \lambda^+$ is an
order-isomorphism then $(\alpha, v') \in
\Pp$.  We write $h[p]=(\alpha, h[v])$.
Moreover, if $q\le p$ then $h[q]\le
h[p]$.
\item 
\label{amalg} {\em (Amalgamation)} For every $p,q
\in 
\Pp$ and 
$\alpha<\lambda^+$, if $p\restr \alpha
\le q$ and $\dom p \intersect \dom q=
\dom p \intersect
\alpha$, 
then there exists $r \in \Pp$ so that
$p,q \le r$.
\ee
\end{defn}

\medskip

It should be remarked that a standard $\lambda^+$-uniform partial
order comes with the additional structure imposed on it by the domain
and restriction functions.  We will call a partial order {\em
$\lambda^+$-uniform} if it is {\em isomorphic} to a standard
$\lambda^+$-uniform partial ordering.  It follows that although a
$\lambda^+$-uniform is isomorphic to a standard one {\em as a partial
order}, there will be an induced notion of domain and restriction. The
elements of such a partial order will be called {\em approximations},
rather than ``conditions'', as we are aiming at a construction in the
ground model.

Observe that $p\restr \alpha=p$ iff
$\dom \,p\includedin
\alpha$. 
 Note also that for $p\le q$ in $\Pp$,
$p\restr
\alpha\le 
q\restr \alpha$.  (As $p\restr
\alpha,q\restr
\alpha\le 
q$, there is $r\le q$ in $\Pp$ with
$p\restr \alpha,q\restr \alpha\le r$ and
$\dom r=\dom p\restr \alpha\union\dom
q\restr
\alpha=\dom 
q\restr \alpha$; hence $r=q\restr
\alpha$ by maximality of $q\restr
\alpha$, and $p \restr \alpha \le
q\restr \alpha$.)

It is important to realize that in intended applications there will be
$\lambda$ many comparable elements of a $\lambda^+$-uniform partial
order which have the same domain (see the first example of the next section).

Typically the only condition that
requires attention in concrete cases is
the amalgamation condition.  
It is therefore useful to have a
weaker version of the amalgamation
property available which is sometimes
more conveniently verified, and which
is equivalent to the full amalgamation
condition in the presence of the other
(trivial) hypotheses.
Such a version is:

\begin{quote}  Weak Amalgamation. {\em
For every
$p,q\in \Pp$, and $\alpha < \lambda^+$, if
$p\res \alpha \le q$, $\dom p
\includedin \alpha+1$,  and $\dom q 
\includedin \alpha$, then there 
exists $r\in \Pp$ with $p,q\le r$.}
\end{quote}

To prove amalgamation from weak amalgamation, we define a continuous
increasing chain of elements $r_\beta \in \Pp$ for $\beta \geq \alpha$
so that
\be \item $\dom(r_\beta) \sub \beta$ and
\item $r_\beta \geq p \res\beta, q \res\beta$.
\ee

Let $r_\alpha = q \res\alpha$.  For limit ordinals, use conditions 3
and 5 of the definition of uniform partial order.

Suppose we have defined $r_\beta$ and $\beta \not \in \dom(p) \cup
\dom(q)$.  Let $r_{\beta +1} = r_\beta$.

If $\beta \in \dom(q) \sm \dom(p)$ then $p \res{\beta +1} = p
\res\beta$.  Apply weak amalgamation to $r_\beta$ and $q\res\beta$.
Using condition 2 of the definition now, we can define $r\res{\beta
+1}$.

If $\beta \in \dom(p) \sm \dom(q)$ then we can apply weak amalgamation
to $p\res{\beta +1}$ and $r_\beta$.

Since these are all the possibilities,  let $\gamma =
\sup(\dom(p) \cup \dom(q))$ and so $r_\gamma \geq p,q$. This
verifies amalgamation.

\setoff 
Notation

For $p,q\in \Pp$ we write $p \lesd q$ to
mean $p \le q$ and $\dom p=\dom q$.
(Here ``sd'' stands for ``same
domain''.)
If $p,q \in \Pp$ then we write $p \perp q$ if $p$ and $q$ are
incompatible i.e. there is no $r$ so that $p \leq r$ and $q \leq r$.

We define the {\em collapse} $p^{\col}$
of an approximation as $h[p]$ where $h$
is the canonical order isomorphism
between $\dom p$ and $\otp(\dom p)$.

\medskip

\noindent{\bf 
Convention}
\smallskip

For the remainder of this section we fix
a standard $\lambda^+$-uniform partial
order $\Pp$, and we let
\[\Pp_\alpha 
= \{ p \in \Pp :
\dom 
p \sub \alpha \}\] for $\alpha <
\lambda^+$.  Note that $\Pp_{\lambda^+} = \Pp$.

\medskip
Be forewarned that the following definition does not follow the
standard set theoretic use of the term ``ideal''.

\begin{defn}
\be \item 
For $\alpha<\lambda^+$, a {\em
$\lambda$-generic ideal} $G$ in $\Pp_\alpha$ is a subset
of $\Pp_\alpha$ satisfying:
\be \item $G$ is closed downward;
\item if $Q \sub G$ and $|Q| < \lambda$ then
$Q$ has an upper bound in $G$; and
\item for every $p \in \Pp_\alpha$, if $p \not
\in G$ then $p$ is incompatible with
some $q \in G$. \ee $\gen(\Pp_\alpha)$
is the set of $\lambda$-generic ideals
of $\Pp_\alpha$.
\item  If $G \in \gen(\Pp_\alpha)$ then
\[\Pp/G = \{ p \in \Pp : p
\mbox{ is compatible with every } r \in G \}.\]
Note that $p \in \Pp/G$ iff $p \restr
\alpha \in G$.
\item We say an increasing sequence $\la g_i : i
< \lambda \ra$ is {\em cofinal} in $G \in
\gen(\Pp_\alpha)$ if $G = \{ r \in
\Pp_\alpha : \mbox{ for some $i$, } r
\le g_i \}$.  Every $G \in
\gen(\Pp_\alpha)$ 
has a cofinal sequence of length
$\lambda$ (possibly constant in
degenerate cases).  We will often write $(g_\delta)_\delta$ to mean
$\la g_\delta : \delta < \lambda\ra$.
\item We will say that $G$ is generic if $G \in \gen(\Pp_\alpha)$ for 
some $\alpha$.
\ee
\end{defn}

\medskip

\begin{lemma}
\label{mingen}
Let $G_i \in \gen(\Pp_{\alpha_i})$ for
$i<\delta$ be an increasing sequence of
sets, and $\alpha=\sup_i \alpha_i$.
Then there is a unique minimal
$\lambda$-generic ideal of
$\Pp_{\alpha}$ containing $\bigcup_{i
<\delta} G_i$. This ideal will be
denoted $G_{<\delta}$.

\end{lemma}

\proof 
We may suppose that $\delta$ is a
regular cardinal, $\delta\le\lambda$.
If $\delta = \lambda$ then it is clear
that $\Union_{i < \delta} G_i \in
\gen(\Pp_\alpha)$.  Suppose now that
$\delta<\lambda$.  For $i<\delta$ fix an
increasing continuous sequence
$(g^i_\gamma)_{\gamma < \lambda}$ cofinal
in $G_i$. Fix $i<j<\delta$. 
There is a club $C_{ij}$ in $\lambda$ such 
that for all $\gamma \in C_{ij}$,
$g^i_\gamma = g^j_\gamma \restr
\alpha_i$.
Let $C = \bigcap_{i < j< \delta}
C_{ij}$.  If $\beta \in C$ then define
$g_\beta = \bigcup_{i < \delta}
g^i_\beta \in \Pp_\alpha$.  Then the
downward closure of $(g_\beta : \beta <
\lambda)$ is the required generic set in
$\Pp_\alpha$. \qed
\bigskip

The notion of $\lambda$-genericity is of
course very weak.  In order to get a
notion adequate for the applications, we
need to formalize the notion of a {\em
uniform} family of dense sets.
\medskip

\begin{defn}\nl
\be \item 
For $\alpha <\lambda^+$ and $G \in
\gen(\Pp_\alpha)$ or $G = \emptyset$ (in which case, in what follows,
read $\Pp$ for $\Pp/G$) we say
\[D : \{(u,w) : u \sub w \in 
\CP_{< \lambda}(\lambda^+) \} \rightarrow
\CP(\Pp) \] is a 
{\em density system}
over $G$ if:
\be \item 
for every $(u,w)$, $D(u,w) \sub \{p
\in 
\Pp/G : \dom p \sub w  \}$,
\item 
for every $p,q \in \Pp/G$, if $p \in
D(u,w)$, $p \le q$ and $\dom q \sub w
$ then $q \in D(u,w)$,
\item 
{\em (Density)} For every $(u,w)$ and
every $p \in \Pp/G$,
with $\dom p \sub w $,
 there is $q\ge p$ in $D(u,w)$; and
\item 
{\em (Uniformity)} For every $(u_1,w_1), (u_2,w_2)$, if $w_1 \cap
\alpha = w_2 \cap \alpha$ and there is an order-isomorphism $h: w_1
\rightarrow w_2$ such that $h[u_1] = u_2$, then for every $p\in
\Pp/G$ with $\dom p\includedin w_1$
\[p 
\in D(u_1,w_1) 
\mbox{ iff } h[p]\in D(u_2,w_2).\] 
\ee

The term ``density system'' will refer
to density systems over
some $G \in
\gen(\Pp_\alpha)$, for some $\alpha$, 
and we write ``0-density system'' for
density system over $\emptyset$.
\medskip

\item 
For $G \in \gen(\Pp_\gamma)$ and $D$ any
density system, we say $G$ {\em meets}
$D$ if for all $u \in \CPll(\gamma)$
there is $v\in \CPll(\gamma)$ so that $u
\sub v$ and $G \intersect D(u,v) \ne
\emptyset$.
\ee
\end{defn}

We give now two examples of density systems which will be important in
the proof of Theorem \ref{main}.  Both examples use the following notion. 
A closed set $X$ of
ordinals will be said to be {\em
$\lambda$-collapsed} if $0 \in X$ and for any
$\alpha\le\sup\,X$,
$[\alpha,\alpha+\lambda]\intersect X\ne
\emptyset$.  An order isomorphism
$h:Y\with X$ between closed sets of
ordinals will be called a {\em
$\lambda$-isometry} if for every pair
$\alpha\le \beta$ in $Y$ and every
$\delta<\lambda$, $\beta=\alpha+\delta$
iff $h(\beta)=h(\alpha)+\delta$.  Every
closed set of ordinals is
$\lambda$-isometric with a unique
$\lambda$-collapsed closed set; the
corresponding $\lambda$-isometry will be
called the $\lambda$-collapse of $Y$,
and more generally the
$\lambda$-collapse of any set $Y$ of
ordinals is defined as the restriction
to $Y$ of the $\lambda$-collapse of its
closure.  Observe that a
$\lambda$-collapsed set of fewer than
$\lambda$ ordinals is bounded below
$\lambda\times
\lambda$ 
(ordinal product).

\begin{example}
\label{ex1}

We shall show that there is a family
$\CD$ of at most $\lambda$ 0-density
systems such that for any $\alpha
<\lambda^+$, if $G\in\gen(\Pp_\alpha)$
meets all $D\in \CD$ then $\Pp/G$ is
again $\lambda^+$-uniform.  (The
amalgamation property must be verified.)
\medskip     

\noindent{\bf Construction}

For $p,q \in \Pp_{\lambda \times \lambda}$ and $\delta <
\lambda
\times 
\lambda$ 
(where $\lambda\times \lambda$ is the
{\em ordinal} product), we define a
density system $D_{p,q,\delta}$ as
follows.  Let $u =(\dom p \union \dom q)
\intersect
\delta$.  For \[u' \sub w' \in
\CPll(\lambda^+),\] if there is an
order-isomorphism $h : w' \rightarrow w
\sub \delta$ with $h[u'] = u$, then let
\[\begin{array}{lcl}
D_{p,q,\delta}(u',w') & = & \{r : \dom r\sub w' 
\mbox{ and either there does
not exist }\\ & & \mbox{\quad } s \ge p,q,h[r] 
\mbox{, or there exists } s\ge p,q \\ & & \mbox{\quad so that } 
s\restr \delta = h[r] \}.
\end{array} 
\]
This definition is independent of the choice of $h$.
     
If there is no such $h$ then let
$D_{p,q,\delta}(u',w') = \{r : \dom p
\sub w' \}$.
We claim that $D_{p,q,\delta}$ is a
0-density system.  It suffices to check
the density condition 
for $u\includedin w\includedin \delta$,
and this is immediate.

\noindent{\bf 
Application}

We will now show that if $G \in
\gen(\Pp_\alpha)$ meets every
density system of the form
$D_{p,q,\delta}$ then $\Pp/G$ is
$\lambda^+$-uniform.  In order to view
$\Pp/G$ as encoded by elements of
$\lambda\times\CPll(\lambda^+)$, we let
$h:\lambda^+\backslash \alpha\with
\lambda^+$ be an order isomorphism, and
replace $(\beta,u)$ in $\Pp/G$ by
$(\beta',h[u\backslash \alpha])$ where
$\beta'$ is just a code for the pair
$(\beta, u\intersect\alpha)$.  We need
only check the amalgamation condition
(\ref{amalg}) of the definition.
     
Let $p,q \in \Pp/G$, $\beta < \lambda^+$
with $p \restr \beta
\le 
q$ and $\dom q \intersect \dom p = \dom
p\intersect \beta$.  We must find $r
\ge p, q$ with $r \in \Pp/G$.  Let $X =
\dom p \union
\dom 
q \union \{ \alpha \}$ and let $h_0:X
\to X'$ be the $\lambda$-collapse of
$X$.  Let $p'= h_0[p], q' = h_0[q]$,
$\alpha' = h_0(\alpha)$, and $u = \dom q
\intersect
\alpha$. 
   Now choose $w \sub \alpha$ with  $|w|<\lambda$ and
$r \in G \intersect D_{p',q',\alpha'}(u,w)$.  Since $X'$ is
$\lambda$-isomorphic with $X$,
we can extend $h_0$ to
an order-isomorphism
\[h 
: X \union w \rightarrow X' \union w'
\mbox{ with } h[w] = w' \sub\alpha'.\]

We claim that there is $s \ge p',q',h[r]$.  It suffices to find some
$s \ge p,q,r$.  Since $p\restr \alpha, q\res\alpha,r$ are all in $G$,
we may take $r'\ge p\restr \alpha, q\res\alpha,r$ in $G$.  Since $q
\in \Pp/G$ and $r' \in G$ then by amalgamation we can find $\hat q \ge q,
r'$ with $\dom \hat q=\dom q\union\dom r'$.  But now $\dom p \intersect
\dom
\,\hat q = \dom p \intersect \beta$ and
$p\restr \beta \leq \hat q$, so we
can find $s\ge p, \hat q$.  This is the desired $s$.
     
As $r \in D_{p',q',\alpha'}(u,w)$, it
now follows that there exists $s \ge p', q'$
so that $s\res\alpha' = h[r]$, and hence
$h^{-1}[s] \ge p,q$ and
$(h^{-1}[s])\restr \alpha = r$.  So
$h^{-1}[s] \in \Pp/G$ and $h^{-1}[s]\ge
p,q$, verifying condition (\ref{amalg})
for $\Pp/G$.
\end{example}
\bigskip

\begin{example}
\label{ex2}
\quad
The next example will be useful in the
following situation.  Suppose we have $G
\in \gen(\Pp_\alpha)$, $\beta>\alpha$,
and we want to build $G' \supseteq G$
with $G' \in \gen(\Pp_\beta)$.  To
ensure the genericity of $G'$ we must
arrange that for all $q \in \Pp_\beta$,
either $q \in G'$ or else $q$ is
incompatible with some $g \in G'$.  We
will find another family of at most
$\lambda\ $ 0-density systems
$D_{p,q,r,\delta}$ which make it
possible to construct a suitable
$G'\containing G$ if $G$ meets all
$D_{p,q,\delta}$ (from Example 1.5)
and $D_{p,q,r,\delta}$.
\medskip

\noindent{\bf Construction}

For $p,q,r \in \Pp_{\lambda \times
\lambda}$, $\delta <
\lambda 
\times
\lambda$ 
such that:
\[p \restr \delta \le r;
\ \dom r\includedin \delta; \mbox{ 
and there does not exist $s \ge p,q,r$,}\] 
we define $D_{p,q,r,\delta} $ as follows 
     
Let $u = (\dom p \union \dom
q)\intersect \delta$.  For $u' \sub w'
\in
\CPll(\lambda^+)$, 
if there is an order-isomorphism $h:w'
\rightarrow 
w$ where $w \sub \delta$ and $h[u'] = u$
then let
\[ 
\begin{array}{ccl}
D_{p,q,r,\delta}(u',w') & = & \{ s :
\mbox{$\dom s \sub w'$ and $h[s]$ is
incompatible with $r$,} \\ & &
\mbox{\quad or $h[s] \ge r$ and there is
some $t \ge p$ so that }\\ & &
\mbox{\quad $t \restr \delta \le h[s]$
and $t$ is incompatible with $q$}\}.
\end{array} 
\]
If there is no such $h$ then let
$D_{p,q,r,\delta}(u',w') = \{ s :
\dom 
s \sub w' \}$.
     
We claim that $D_{p,q,r,\delta}$ is a
0-density system.  Again we check only
the density condition for $u\includedin w 
\includedin \delta$.  
So suppose we
have $s \in \Pp$, $\dom s \sub w$, and
$s$ is compatible with $r$. We seek
$s'\ge s$ in $D_{p,q,r,\delta}(u,w)$.
     
Choose $s'\ge r,s$ with domain $\dom r\union\dom s$; 
so $\dom s'\includedin \delta$.  
Then $s' \ge r \ge p \restr \delta$ and 
$\dom s' \intersect\dom p = \dom p \intersect \delta$, 
so we can choose $t \ge s', p$
so that $\dom(t) = \dom s'\union\dom p$, 
and hence $t$ is incompatible with $q$
(since there is no $t' \ge p,q,r$).  
Now $t \restr \delta \ge s' \ge r,s$, so
if $s''= t \restr \delta$ then
$s'' \ge r,s$, and 
$s''\in D_{p,q,r,\delta}(u',w')$ as desired.

\bigspace
\end{example}
\bigskip
\noindent{\bf 
Application}
\medskip

We return to the situation in which we
have $G \in
\gen(\Pp_\alpha)$, 
$\beta>\alpha$, and we want to build $G'
\supseteq 
G$ with $G' \in \gen(\Pp_\beta)$,
assuming that $G$ meets all
$D_{p,q,\delta}$ and $D_{p,q,r,\delta}$.
We will naturally take $G'$ to be the
downward closure of a sequence
$(g_i)_{i<\lambda}$ which is constructed
inductively, taking suprema at limit
ordinals.  At successor stages, suppose
that the $i$-th term of our sequence has
just been constructed, and let $p=g_i$.
Suppose $q \in \Pp_\beta$ is fixed.
We wish to ``decide'' $q$: that is, we
seek $\hat p \ge p$ so that either $\hat
p$ is incompatible with $q$, or else
$\hat p \ge q$.

If $p$ is already incompatible with $q$
then let $\hat p=p$.  Otherwise, let $X
= \dom p \union \dom q
\union 
\{ \alpha\}$ and let 
$h : X \rightarrow X'$ be the
$\lambda$-collapse of $X$.  Let $p' =
h[p]$, $q' = h[q]$, and $\alpha' =
h(\alpha)$.  If $u = X \intersect
\alpha$, choose $w \supseteq u$ and $r
\in G \intersect
D_{p',q',\alpha'}(u,w)$.  Extend $h$ to
an order-preserving function from $X
\union w$ to $X' \union w'\includedin
\alpha'$, and let $r'=h[r]$.
     
Suppose first that there is some $s \ge
p',q'$ with $s \restr \alpha' \le r'$.
We may suppose that $\dom s = \dom
p'\union\dom q'$.  In this case let
$\hat p = h^{-1}[s]$.  As $\hat p \restr
\alpha \le r$, we have $\hat p \in
\Pp/G$, and $q$ is decided by $\hat p$.
     
Now suppose alternatively that there is
no $s \ge p',q',r'$.  We may assume that
$p \restr \alpha \le r$ since $p\restr
\alpha\in G$ and $G$ is directed.  Let:
\[Y = \dom p \union \dom q \union \dom r
\union \{ \alpha \},\] and let 
$k: Y\rightarrow Y''$ be the
$\lambda$-collapse of $Y$.  Let $ 
p'' = k(p), q'' = k(q),r'' =k(r)$, and
$\alpha'' = k(\alpha)$.  Then $p'' \restr
{\alpha''} \le r''$, and there is no $s
\ge p'', q'', r''$.

\begin{sloppypar}     
Let $v = (\dom q \union \dom r)
\intersect \alpha$, and choose $z
\supseteq v$ and $s \in G \intersect D_{
p'', q'', r'', \alpha''}(v,z)$.  We can
extend $k$ to an order-isomorphism from
$Y \union z$ to $Y'' \union z''$ with
$k[z] = z'' \sub {\alpha''}$. Let
$s''=k[s]$.
\end{sloppypar}
     
Certainly $r''$ and $s''$ are compatible
since $r,s \in G$.  As $s$ belongs to $D_{ p'',
q'', r'', \alpha''}(u,z)$, we have $k[s]
\ge r''$, and there is some $t'' \ge
p''$ so that $t'' \res {\alpha''} \le
s''$ and $t''$ is incompatible with $
q''$; in other words, $s \ge r$, and
there is some $\hat p \ge p$ so that
$\hat p \restr \alpha
\le 
s$ and $\hat p$ is incompatible with
$q$. Then $\hat p \in \Pp/G$, and $\hat
p$ decides $q$.\qed
\bigskip

We now introduce the {\em genericity
game}. Our main theorem will state that
the second player has a winning strategy
in this game, under certain
combinatorial conditions.

\begin{defn} 
Let $\Pp$ be a $\lambda^+$-uniform
partial order.  The {\em genericity
game} for $\Pp$ is the two-player game
of length $\lambda^+$ played according to
the following rules:
\be\item 
At the $\alpha^\th$ move, player II will
have chosen an increasing sequence of
ordinals $\zeta_\beta<\lambda^+$, and
will have defined an increasing sequence
of $\lambda$-generic ideals $G_\beta$ on
$\Pp_{\zeta_\beta}$ for all $\beta <
\alpha$. 
Player I will choose an element
$g_\alpha \in
\Pp/G_{<\alpha}$ 
and will also choose at most $\lambda$
density systems $D^\alpha_i$ over
$G_{<\alpha}$.  Note that $G_{<\alpha} \in
\gen(\Pp_{\zeta_{<\alpha}})$ by Lemma \ref{mingen}.

\item 
After player I has played his
$\alpha^\th$ move, player II will pick
an ordinal $\zeta_\alpha$ and a
$\lambda$-generic ideal of $\Pp_{\zeta_\alpha}$.
\ee

Player II wins the $\Pp$-game if the
sequences $\zeta_\alpha$ and $G_\alpha$
are increasing, and for all $\alpha$,
and all indices $i$ occurring at stage
$\alpha$: $g_\alpha \in G_\alpha$,
and for all $\beta \ge \alpha$, $G_\beta$ meets
$D^\alpha_i$.
\end{defn}

Our main theorem uses the following
combinatorial principle.
     
\begin{defn} 
\label{defdl} Suppose $\lam$ is a regular cardinal.
$\Dl_\lambda$ asserts that there are sets
$\CA_\alpha \sub
\CP(\alpha)$, 
$|\CA_\alpha| < \lambda$ for every
$\alpha < \lambda$, such that for all $A
\sub \lambda$:
\[ 
\{ \alpha \in \lambda : A \intersect \alpha \in \CA_\alpha\} \mbox{ is
stationary.}\]
\end{defn}

Easily, $\Diamond_\lambda$ or $\lambda$ strongly inaccessible (or even
$\lambda = \aleph_0$) implies $\Dl_\lambda$.  Also, Kunen showed that
$\Dl_{\lambda^+}$ implies $\Diamond_{\lambda^+}$.  Gregory has shown
that if GCH holds and $\cof(\kappa) > \aleph_0$ then
$\Diamond_{\kappa^+}$ holds. 
It is useful to note that $\Dl_\lambda$ implies $\lambda^{<\lambda} =
\lambda$.

\begin{theorem}
\label{main}
$\Dl_\lambda$ implies that player II has
a winning strategy for the $\Pp$-game.
\end{theorem}

This theorem will be proved in \S\S3-5.
We illustrate its use in the next
section.
\page

\section{Illustrative application}

In this section we give an example of an application of
the combinatorial principle described in section 1.



In \cite{MagMal}, Magidor and
Malitz introduce a logic $\Cl^{<\omega}$ which
has a new quantifier $Q^n$ for each $n \in
\omega$, in
addition to the usual first order connectives and
quantifiers.
The $\kappa$-interpretation of the formula
$Q^n \varphi (x_1,\ldots,x_n,\bar y)$ is
\begin{quote}
``there is a set $A$ of cardinality $\kappa$ so that
for any  
$x_1, \ldots, x_n \in A$, $\varphi(x_1, \ldots, x_n,\bar y)$ holds.''
\end{quote}
They then give a list of axioms which are sound for the
$\kappa$-interpretation when $\kappa$ is regular, and
show that these axioms are complete for the
$\aleph_1$-interpretation under the assumption of
$\Diamond_{\aleph_1}$.  They ask  whether 
these axioms are
complete for the $\lambda^+$-interpretation.  We 
will show that their axioms are  complete when both
$\Dl_\lambda$ and $\Diamond_{\lambda^+}(\{
\delta < \lambda^+ : \mbox{ cf}(\delta) = \lambda \})$
hold. This will explain a remark at the end of
\cite{secondIII}. See Hodges (\cite{Hodges}) for a treatment in the
same vein for the $\aleph_1$-interpretation.

Fix a complete $\Cl^{<\omega}$ theory $T$, $|T| \leq
\lambda$.  Let $Q = Q^1$.  We may assume that
that associated to each formula
$\varphi$ with free variables $x_1, \ldots, x_n, y_1,
\ldots, y_m$, $\Cl$ contains an $(m+1)$-ary function
$F_\varphi$, so that $T$ proves
$Qx(x=x)$ and, 
for any fixed $y_1,
\ldots y_m$, $F_\varphi(-,y_1, \ldots, y_m)$ is one-to-one
and
\begin{eqnarray*}
\lefteqn{ Q^n\bar x\varphi(x_1, \ldots, x_n, y_1, \ldots, y_m) \rightarrow} \\
& & \forall z_1 \ldots z_n (\bigwedge_{i<j} z_i \neq z_j \rightarrow
\bigwedge_{\sigma\in Sym(n)}\varphi(
F_\varphi( z_{\sigma(1)},\bar y), \ldots, F_\varphi(z_{\sigma(n)},\bar y),\bar
y)). \end{eqnarray*}
Strictly, it is not necessary to make this conservative extension to
our language and theory but it is convenient when handling the
inductive step corresponding to $Q^n$.

We add new constants $\{ y_\alpha :\alpha < 
\lambda^+ \}$ and
$\{x^\alpha_i : \alpha < \lambda^+, i < \lambda \}$
to $\Cl$, to obtain a language
$\Cl_1$.  
The set of constants $\{y_\alpha\} \cup \{x^\alpha_i:i <
\lambda\}$ is called the set of $\alpha$-constants and
$y_\alpha$ is called the {\em special} $\alpha$-constant.  
A constant is
said to be a $w$-constant if it is a $\beta$-constant for some
$\beta \in w$; in particular a constant is a
$(<\alpha)$-constant 
if it is a $\beta$-constant for some $\beta
< \alpha$.

We define a partial order $\Bbb P$ as follows: $p \in \Bbb P$
iff \be \item $p$ is a set of $\Cl_1^{<\omega}$ sentences
consistent with $T$;
\item $|p| < \lambda$;
\item $p$ is closed under conjunction and existential
quantification; and
\item if $\varphi(y_\alpha,\bar z) \in p$ and the 
$\bar z$ are 
($<\alpha$)-constants, then $Qy \varphi(y, \bar z) \in p$. \ee

We now indicate how
$\Bbb P$ may be viewed 
as a standard $\lambda^+$-uniform partial order.
We order $\Bbb P$ by inclusion.  
Let $\Pp_\alpha$ be 
\[\{p\in \Pp: \hbox{all constants occurring in a
formula of $p$ are $(<\alpha)$-constants}\}.\]
The elements of $\Pp_\lambda$ will be called {\em templates}.
For any template $p$, there is a least $\beta$ so that all
formulas in $p$ use only constants from $\{y_i : i < \beta \}
\cup \{x^i_j : i < \beta, j < \lambda \}$.  Call this
$\beta_p$.

For any template $p$ and any $w \sub \lambda^+$ so that
$\otp(w) = \beta_p$, fix an order-isomorphism $h:\beta_p
\rightarrow w$.  Define $p(w)$ as the set of formulas obtained
by replacing $x^i_j$ and $y_i$ by $x^{h(i)}_j$ and $y_{h(i)}$
respectively for $i < \beta_p$.  
Every element of $\Bbb P$ can be obtained in this way from
a template.

Let $\iota$ be any bijection between the set of templates and
$\lambda$.  Identify $\Pp$ with the set $\{(\iota[p],w): 
p\in \Pp_\lambda, w\in \CPll(\lambda^+) \mbox{ where } \otp(w) =
\beta_p\}$ by sending $(\iota[p],w)$ to $p(w)$.  Throughout the rest
of this section we will treat $\Bbb P$ as if it were in standard for
although in practice we will use its original definition.
We claim that $\Bbb P$ is $\lambda^+$-uniform; 
it suffices to check the amalgamation condition \ref{amalg}.

The following notation will be convenient.
If $\varphi(y_{\alpha_1}, \bar z_1,y_{\alpha_2}, \bar z_2,
\ldots,y_{\alpha_n}, \bar z_n )$ is a formula with $\alpha_1 > \alpha_2 >
\ldots > \alpha_n$, and $\bar z_i$ is a collection of
$[\alpha_{i+1},\alpha_i)$-constants, then the string $S$ of
quantifiers:
\[ \exists \bar x_n Q y_n \ldots \exists \bar x_1 Q y_1 \] is
called {\em standard} for $\varphi$ where the $x$'s quantify over the
$z$'s.  Its dual is denoted $S^*$: 
\[ \forall \bar x_n \neg Q y_n \neg \ldots \forall \bar x_1
\neg Q y_1 \neg \] 
If $p$ is a set of fewer than $\lambda$ formulas
of $\Cl_1^{<\omega}$ which is closed under conjunction,
then the following are equivalent:
\be
\item $p
\sub q$ for some $q \in \Bbb P$;
\item $S\varphi \in T$ for all
$\varphi \in p$ where $S$ is standard for $\varphi$.
\ee

For $p \in \Bbb P$ and $\alpha < \lambda^+$, we have:
\[ p \restriction \alpha = \{ \varphi \in p : \mbox{ all
constants in } \varphi \mbox{ are $<\alpha$-constants } \}. \] To show
that $\Pp$ satisfies amalgamation, we will show that it satisfies weak
amalgamation.  Suppose $p \in \Bbb P_{\alpha +1}$, $q \in \Pp_\alpha$
and $p\res\alpha \leq q$.

Suppose $\varphi(\bar x,y_\alpha,\bar z) \in p$ where $\bar x$ is all the $\alpha$-variables except $y_\alpha$ and $\bar z$ is the $<\alpha$-variables.  Then
\[ Qy\exists \bar x \varphi(\bar x,y,\bar z) \in p\res_\alpha .\]
If $\psi \in q$ then 
\[SQy\exists \bar x(\psi \wedge \varphi) \] where $S$ is a standard
sequence, is equivalent to 
\[S(\psi \wedge Qy\exists \bar x \varphi) .\]  Since both of the
conjuncts are in $q$, this last sentence is in $T$.  This verifies
weak amalgamation.

Now the strategy is to build a  set $G$ which is the union of generics
so that the constant structure derived from $G$
will form a model of $T$ under the
$\lambda^+$-interpretation. More precisely, we introduce an
equivalence relation $\sim$ on the set of nonspecial constants 
$A = \{ x^\alpha_j : \alpha < \lambda^+, j < \lambda \}$ by:
\[a \sim b \mbox{ iff ``$a=b$''} \in G.\]
Let $\bar G = \{ a/\!\!\sim : a \in A \}$ and define the
functions and relations on $\bar G$ in the usual manner.  We
want to ensure that for any formula
$\varphi$ in $\Cl_1^{<\omega}$ we will have:
\begin{equation} \label{equation1}
\bar G \models \varphi(a_1/\!\!\sim , \ldots ,
a_n/\!\!\sim) \mbox{ iff } \varphi(a_1,\ldots,a_n) \in G.
\end{equation}
If (1) is true, its proof naturally proceeds by induction on the 
complexity of formulas.  We now describe a strategy for Player
I in the genericity game which can only be defeated by
achieving (1).  In other words, we will specify density systems
and elements $g\in \Pp$, to be played by Player I, such that
a proper response by Player II ensures that $G$ allows an 
inductive
argument for (1) to be carried out.  Our
discussion will be somewhat informal, stopping well 
short of actually writing down the density systems in
many cases.

We begin with the treatment of the ordinary
existential quantifier.
Whenever $\exists x \varphi(x,\bar z) \in G$ we will want 
(eventually) to have some
$a \in A$ so that $\varphi(a,\bar z) \in G$.  In particular,
for every $\alpha$ there will be some $a \in A$ so
that $y_\alpha = a \in G$.  The density systems which ensure
this condition is met will in fact be
0-density systems.

Next  we consider the quantifier $Q$. 
For each formula
$Qx \varphi(x)$
which is put into $G$, at cofinally many 
subsequent stages we wish to add the
formula $\varphi(y)$ for an unused special constant $y$.
The first player will play such formulas as ``$g_\alpha$''
from time to time.  We will also have to deal with the case
in which $\neg Qx\varphi(x)$, and we will return
to this in a moment.

We now consider the quantifier $Q^n$.
Suppose that the formula $Q^n \bar x \varphi(\bar x,\bar y)$ is
in $G$ at some stage.
This is where we use the function $F_\varphi$.  If $\bar G$ is a model
of $T$ then it has cardinality $\lambda^+$.  Moreover,
$F_\varphi(-,\bar y)$ is one-to-one.  Since $Q^n \bar x \varphi(\bar
x,\bar y)$ is in $G$, so is
\[ \forall z_1 \ldots z_n (\bigwedge_{i<j} z_i \neq z_j \rightarrow
\bigwedge_{\sigma\in Sym(n)}\varphi(
F_\varphi( z_{\sigma(1)},\bar y), \ldots, F_\varphi(z_{\sigma(n)},\bar y),\bar
y)).\]
It follows that the range of $F_\varphi(-,\bar y)$ is homogeneous for
$\varphi$.

We are now left with the cases in which formulas
of 
the form $\neg Q^n x \varphi$ ($n\ge 1$) are placed
in $G$.  We deal first with the case $n=1$.  For this case,
we define
a number of density systems depending on the following parameters:

\be \item $j, j_0, \ldots, j_{m-1} < \lambda$;
\item a formula $\varphi(x,y_0,\ldots,y_{m-1})$; and
\item a function $f:m \rightarrow m$.
\ee

We associate with these data a density system $D$. If
$\otp(u) \neq m+1$, we let $D(u,w)$ be degenerate:
\[D(u,w)= \{ p \in \Bbb P : \dom(p) \sub w \}.\]
If $\otp(u) = m+1$ then let $g : m+1 \rightarrow u$ be an order
preserving map, let $h = gf$ and set $\beta = g(m)=\max u$, and:
\[ \psi(x) =
\varphi(x,x^{h(0)}_{j_0},\ldots,x^{h(m-1)}_{j_{m-1}}).\]
We will then let $D(u,w)$ consist of those  $p \in \Bbb P$
for which, setting $\alpha=\min(w)$, we have:
\be\item$\dom(p) \sub w$;
\item\mbox If $\neg Qx \psi(x) \in p$, then either
$\neg \psi(x^\beta_j) \in p$ or $x^\beta_j = x^\alpha_i \in p$
for some $i < \lambda$.
\ee

\begin{sloppypar}
We shall verify the density condition on $D$.
Suppose $q \in \Bbb P$ and $\neg Qx
\psi(x) \in q$.  The extension of $q$ we are about to
construct will only involve the adjunction of formulas with 
($\leq\beta$)-constants, so we may assume that $q$ itself
contains only ($\leq \beta$)-constants. 
\end{sloppypar}

If we
cannot complete $q \cup \{ \neg
\psi(x^\beta_j) \}$ to an element of $\Bbb P$,
then there is some $\chi \in q$ so that:
\[ S \exists x ( \chi \wedge \neg \psi(x) ) \not \in T \]
where $S \exists x$ is a standard sequence for the
formula $\chi\wedge\neg\psi$.  Note that by the assumption that $\beta$ is
the maximal element of $\dom(q)$, we may assume that
the final quantifier in the standard sequence is an
existential quantifier on the constant $x$ in $\psi$.

By the axioms for the $Q$-quantifier, for any $\theta \in q$ 
such that $T\proves \theta\implies \chi$,
\[ S \exists x ( \theta \wedge \psi(x)) \in T. \]  
As $\neg Q x \psi(x) \in q$, repeated use of 
the $Q$-quantifier axiom:
\[ Qx \exists y \Delta(x,y) \rightarrow \exists x Qy
\Delta(x,y) \vee Qy \exists x \Delta(x,y) \] 
shows that
$\exists x S (\theta \wedge \psi(x)) \in
T$.

If we now choose a constant $x^\alpha_i$ not occurring
in $q$, 
where $\alpha = \min(w)$, one can conclude that $q \cup
\{ x^\beta_j = x^\alpha_i \}$ can be completed to an
element of $\Bbb P$.

It is easy to see that if the foregoing density systems are met, then
we can carry out the argument from right to left in condition (1)
above for $\varphi = Qx \psi(x)$.  We turn now to the treatment of the quantifiers $Q^n$ for
$n>1$.

By applying  Fodor's Lemma to the map sending $\delta$ to
$\dom(f(\delta)) \cap
\delta$ we obtain:
\bl \label{Fodor} If $S \sub \{ \delta : \cof (\delta) = \lambda \}$  is stationary
and $f:S \rightarrow \Bbb P$ then there is a stationary $S' \sub S$, a
template $p$ and $\sigma < \lambda^+$ so that for $\delta \in S'$,
$f(\delta) = p(w_\delta)$ where $w_\delta = \dom(f(\delta))$ and
$w_\delta \cap \delta \sub \sigma$. \el

It will be convenient to treat conditions as if they were single
formulas.  
Extending our previous notation,
for $p \in \Bbb P$ and $S$ a standard sequence 
covering some of the variables in
$p$, we will write $S(p)$
for the set:
\[ \{ S_\varphi \varphi : \varphi \in p \} \] 
where $S_\varphi$ is the standard sequence for $\varphi$ which we
think of as a subsequence of the possibly infinite standard sequence
$S$.

Let $(A_\delta)_{\cof(\delta) = \lambda}$ be a $\Diamond$-sequence.
For $u,v$ sets of ordinals, we write $u < v$ if for all 
$\beta \in u$, $\beta <\min v$.

If $\cof(\delta) = \lambda$ and $G_\delta \in \gen(\Bbb P_\delta)$,
we will define certain associated density systems over $G_\delta$ which
depend on the 
following additional parameters:
\be \item an $i < \lambda$;
\item a formula $\varphi(x_1,\ldots,x_n,\bar y)$ (we will suppress the
$\bar y$);
\item some $k$ with  $0 \leq k < n$;
\item templates $p_1,\ldots,p_k$; and
\item ordinals $\gamma_j < \beta_{p_j}$ for $1 \leq j \leq k$. \ee

The density system $D$ that depends on this particular set
of parameters will be taken to have $D(u,w)$ degenerate
unless:
\be \item $u = \{ \zeta \} \cup \displaystyle{\bigcup_{1 \leq j \leq k} w_j}$;
\item $\delta < \zeta < w_k \sm \delta < \ldots < w_1 \sm \delta$;
\item $w_j \cong \beta_{p_j}$; 
\item $\displaystyle{\bigcup_{1\leq j \leq k } p_j(w_j)}$ can be extended to a member
of $\Pp$;
\ee
in which case we adopt the following notation. Let
$\zeta_j$ be the $\gamma_j^\th$ element of
$w_j$, and write $z^\beta$ for $x^\beta_i$.
Note that since $\gamma_j < \beta_{p_j}$ we will have $\zeta_j >
\delta$ and hence $\zeta < \zeta_k < \ldots < \zeta_1$.
Define the set $r(\alpha_1,\ldots,\alpha_{n-k-1})$ for $\alpha_1 <
\ldots < \alpha_{n-k-1} \in A_\delta$ to be
\[S(\bigwedge_{1 \leq j \leq k} p_j(w_j) \wedge \neg
\varphi(z^{\alpha_1},\ldots,z^{\alpha_{n-k-1}},z^\zeta,z^{\zeta_k},\ldots,z^{\zeta_1}))\]
where $S$ covers all the ($>\zeta$)-variables.

We now define $D(u,w)$ as the set of $q\in \Pp/G$ with
$\dom(q) \sub w$ which satisfy one of the following three
conditions:

\be 
\item $q \perp \displaystyle{\bigcup_{1 \leq j \leq k} p_j(w_j)}$; or
\item $ \displaystyle{\bigcup_{1 \leq j \leq k} p_j(w_j)} \sub q$ and for some $\alpha_1 <
\ldots < \alpha_{n-k-1} \in A_\delta$, \\
$r(\alpha_1,\ldots,\alpha_{n-k-1}) \sub q$; or
\item $ \displaystyle{\bigcup_{1 \leq j \leq k} p_j(w_j)} \sub q$ and for all $\alpha_1 <
\ldots < \alpha_{n-k-1} \in A_\delta$:
\begin{eqnarray*}
\lefteqn{ S^*_0(q \rightarrow S_k^*(p_k(w_k) \rightarrow \ldots
S^*_1(p_1(w_1) \rightarrow }\\
& & \varphi(z^{\alpha_1},\ldots,z^{\alpha_{n-k-1}},z^\zeta,z^{\zeta_k},\ldots,z^{\zeta_1}))\ldots))
\in G_\delta \} \end{eqnarray*} 
\ee

The third condition means that for every $\alpha_1 <
\ldots < \alpha_{n-k-1} \in A_\delta$, there is a $\chi \in q$ and
$\psi_j \in p_j(w_j)$ so that
\[ S^*_0(\chi \rightarrow S_k^*(\psi_k \rightarrow \ldots
S^*_1(\psi_1 \rightarrow
\varphi(z^{\alpha_1},\ldots,z^{\alpha_{n-k-1}},z^\zeta,z^{\zeta_k},\ldots,z^{\zeta_1}))\ldots))
\in G_\delta  \] where $S_j$ covers all the $(\geq
\delta)$-variables in $\psi_j$ for $j > 0$, and $S_0$ covers all the
$(\geq \delta)$-variables in $\chi$.  Notice that the only overlap among
the variables occur in the $(<\delta)$-variables.

Now suppose $\bar G \models Q^n \bar x \varphi(\bar x, \bar
a/\!\!\sim)$.  We would like to argue that $Q^n\bar x \varphi(\bar x,
\bar a) \in G$.  
For convenience we will suppress the parameters $\bar a$.
We may also assume that $T \vdash
\varphi(\bar x) \rightarrow \bigwedge_{i<j} x_i \neq x_j$.

Since $\bar G \models Q^n \bar x \varphi(\bar x)$, there is  a
$\lambda^+$-homogeneous subset $B\includedin \bar G$ for $\varphi$.  
We may assume there is an $i < \lambda$ so that
every $b \in B$ is of the form $x^\alpha_i/\!\!\sim$ for some
$\alpha$.  Let $A = \{ \alpha : x^\alpha_i/\!\!\sim \in B$ and $\alpha$ is the
least such in a given $\sim$-class $\}$.
For any $\delta$ so that $\cof(\delta) = \lambda$, let $\zeta_\delta =
\min(A \sm A_\delta)$.  Note that if $A \cap \delta = A_\delta$ then
$\zeta_\delta > \delta$.

We will now produce the following data.  There will be:
\be
\item stationary sets
$S_k$ for $0 \leq k \leq n$ with $S_{k+1} \sub S_k$ for all $k < n$
and \[S_0 = \{ \delta : \cof(\delta) = \lambda, \mbox{ and } A \cap
\delta = A_\delta \};\]
\item templates $p_k$ for $1 \leq k \leq n$, and ordinals
$\gamma_k$ so that $\gamma_k < \beta_{p_k}$;
\item
a domain $w^k_\delta$ of the same
order type as $\beta_{p_k}$ for each $\delta \in S_k$; 
a $\sigma_k < S_k$ so that if $\delta
\in S_k$ then $w^k_\delta \cap \delta \sub \sigma_k$;
let $\zeta^k_\delta$ be the $\gamma_k^\th$
element of $w^k_\delta$;
\item For $0 \leq k < n$, if
$\delta < \delta_k < \ldots < \delta_1 \in S_k$ are chosen so that
\[\zeta_\delta^{k+1} < w^k_{\delta_k} \sm \delta_k < \ldots < w^1_{\delta_1}
\sm \delta_1\] and $D$ is the density system over $G_\delta$
corresponding to $i,
\varphi, k, p_1, \ldots, p_k$ and $\gamma_1,\ldots,\gamma_k$, then 
\[p_{k+1}(w^{k+1}_\delta) \in D( 
\{\zeta^{k+1}_\delta \} \cup 
w^k_{\delta_k} \cup \ldots \cup w^1_{\delta_1},w^{k+1}_\delta) \cap
G\]
\ee

Using lemma \ref{Fodor} and the fact that $G$ meets all the density
systems introduced at stages $\delta \in S_0$, this is straightforward.

Now suppose $\delta_n < \ldots < \delta_1 \in S_n$, so that
$w^n_{\delta_n} \sm \delta_n < \ldots < w^1_{\delta_1} \sm \delta_1$.
Let $q_k = p_k(w^k_{\delta_k})$.

Since $B$ is a homogeneous set for $\varphi$, it follows that
$\varphi(z^{\alpha_{\delta_n}},\ldots,z^{\alpha_2},z^{\zeta_{\delta_1}})
\in  G$, for every $\alpha_n < \ldots < \alpha_2 < \zeta_{\delta_1}$
in $A_\delta$. Since $q_1 \in G$, using the density systems defined
before, we conclude that
\[S^*_1(q_1 \rightarrow
\varphi(z^{\zeta_{\delta_n}},\ldots,z^{\zeta_{\delta_1}})) \in G\]
where $S_1$ covers the $w^1_{\delta_1}$-variables.  Proceeding by
induction and using the definition of the density systems, we conclude that
\[S^*_n(q_n \rightarrow S^*_{n-1}(q_{n-1} \rightarrow \ldots S^*_1(q_1
\rightarrow \varphi(z^{\zeta_{\delta_n}},\ldots,z^{\zeta_{\delta_1}})))) \in G\]
where $S_j$ covers the $w^j_{\delta_j}$-variables.  

Of course, $S_n q_n \sub G$, so by the Magidor-Malitz axioms,
$Q^n \bar x\varphi(\bar x) \in G$ and we finish.



\page
\section{Commitments}

In this section we begin the proof of
Theorem \ref{main}.  Our main goal at
present is to formulate a precise notion
of a ``commitment'' (that is, a
commitment to enter a dense set -- or in
model theoretic terms, to omit a type).
We will also formulate the main
properties of these commitments, to be
proved in \S\S4-5, and we show how to
derive Theorem \ref{main} from these
facts.

Before getting into the details, we give
an outline of the proof of Theorem
\ref{main}.
\bigskip

\noindent{\bf General overview}
\medskip

Suppose that we wish to meet only the following very simple
constraints.  We are given some $0$-density systems $D_i$ over
for $i<\lambda$, and some $g_0 \in \Bbb P$, and we seek a
$\lambda$-generic ideal $G_0$ containing $g_0$, and 
meeting each $D_i$.
Let $\beta = \lambda \union \sup(\dom(g_0))$, and enumerate
$\Pp_{\beta}$ as $(r_i : i < \lambda)$.  Then we may construct $G_0$
by generating an increasing sequence $(g_\delta)_{\delta<\lambda}$
beginning with the specified $g_0$, and taking $G_0$ to be the
downward closure of $(g_\delta)$.  We will run through this in some
detail.

Our first obligation is to make $G_0$
$\lambda$-generic in $\Pp_{\beta}$.
We will say that $r \in \Pp_{\beta}$
has been {\em decided} if we have chosen
some $g_\delta\in\Pp_\beta$ so that
either $r \perp g_\delta$ or else $r
\leq g_\delta$.  If the sequence
$(g_\delta)_{\delta<\lambda}$ ultimately
decides every $r\in
\Pp_{\beta}$, 
then $G_0$ will be $\lambda$-generic in
$\Pp_{\beta}$.  At stage $\delta+1$ we
will ensure that $r_\delta$ is decided.
This takes care of the basic
$\lambda$-genericity requirement.  At
limit stages we can let $g_\delta$ be
anything greater than $g_{<\delta}$. We
will also take pains at limit stages to
meet the specified density systems
$D_i$.  We enumerate the pairs $(u,
D_i)$ with $u\in
\CP_{<\lambda}(\beta)$, 
using $\lambda^{<\lambda}=\lambda$,
assigning one such pair to each limit
ordinal $\delta<\lambda$.  Suppose that
$(u,D_i)$ is considered at stage
$\delta$.  Let $v=\dom g_{<\delta}$.  By
the density condition on $D_i$, we can
find $g_\delta\ge g_{<\delta}$ with
$g_\delta\in D_i(u,v)$.
     
Thus it is easy to deal with $\lambda$ constraints of the type
arising in one play of our genericity game.  Our strategy in that
game will rely on this sort of straightforward ``do what you must
when you have the time'' approach, but will encounter
difficulties in ``keeping up'' at limit stages in the game.  We
will use $\Dl_\lambda$ to ``guess'' what additional commitments
should be made with regard to various density systems $D_i$, so
that any generic set which we construct subsequently which meets
these commitments will meet each $D_i$.  The commitments
themselves retain the feature that each of them can easily be met
when necessary; deciding when these commitments should be met
requires another use of $\Dl_\lambda$.
     
At stage 0, player I selects some
density systems, to which we may add all
the density systems from examples
\ref{ex1} 
and \ref{ex2}.  From these
we construct some stage 0 commitments
$\up0\pp$, and a $G_0$ meeting
$\up0\pp$.

At stage $\delta$ in the play of the
game, Player II is attempting to extend
$G_{<\delta}$ to a suitable $G_\delta$.
(At limit stages we also will need to
check that $G_\ld$ continues to meet
suitable commitments).  Since $G_\ld$
meets all the previous commitments, in
particular it meets all the density
systems of examples \ref{ex1} and
\ref{ex2}, 
and therefore $\Pp/G_{<\delta}$ is
$\lambda^+$-uniform.  Consequently the
construction of $G_0$ described at the
outset also works in $\Pp/G_{<\delta}$.
Hence we need only construct new
commitments $\up\delta\pp$, add them to
our previous commitments, and construct
$G_\delta$ meeting $\up\delta\pp$ as
above.  In this way, Player II wins the
game.
     
There is a certain difficulty involved
in coping with the freedom enjoyed by
Player I (in terms of obligations
accumulating at limit stages in the
game). There are a priori $\lambda^+$ sets
$u \in
\CP_{<\lambda}(\lambda^+)$ that 
may require attention. On the other
hand, at a given stage $\delta$ we are
only prepared to consider fewer than
$\lambda$ such sets.  However, by
uniformity, it will be sufficient to
consider pairs $(u,w) \in
\CP_{<\lambda}(\beta + 
\lambda) \times
\CP_{<\lambda}(\beta + 
\lambda)$,  and hence $\lambda$ such pairs suffice.
This still leaves Player II at a
disadvantage, but with the aid of
$\Dl_\lambda$, at limit stages we will
guess a relevant set of $u$'s of size
less than $\lambda$.
     
It remains to show that this strategy
can be implemented, and works.

We introduce the notion of basic data which will be provided by
$\Dl_\lambda$. 

\bd \label{defbasic}  A {\em collection of basic data} will contain
\be \item trees $T_\delta$,
subsets of $\Pp_\lambda$  (but not suborders), with orders $<_\delta$,
for every $\delta < \lambda$; 
\item for every generic set $G \in \gen(\Pp_\alpha)$ for some $\alpha
< \lam^+$, two stationary subsets of $\lambda$,
$S(G)$ and $S'(G)$ and a club $C$ so that $C \cap S'(G) \sub S(G)$; and
\item for every $\delta < \lam$, a set $U_\delta \sub \CPll(\lam)$
\ee
with the following properties
\be \item $|T_\delta| < \lam$, $|U_\delta| < \lam$ for every $\delta <
\lam$,
\item if $p \in T_\delta$ then $\len(p) = \alpha$,
\item if $p \leq_\delta q$ and $\len(p) = \alpha$ then $p = q
\res \alpha$,
\item if $p \in T_\delta$ and $\alpha \leq \dom(p)$ then $p\res\alpha
\in T_\delta$,
\item if $(g_\delta)_\delta$ is a cofinal sequence for a generic set
$G \in \gen(\Pp_\alpha)$ then there is a club $C$ so that for $\delta \in C \cap S(G)$,
$(g_\ld)^\col \in T_\delta$,
\item \label{remark} if $G$ and $G'$ are generic sets so that $G \sub
G'$ then there is a club $C$ so that $C \cap S(G') \sub S(G)$ and
\item \label{universality} (oracle property) 
for $\alpha < \lam^+$ and $G \in \gen(\Pp_\alpha)$, $u \in
\CP_\ll(\alpha)$ and $\alpha = \bigcup_{\delta < \lam} w_\delta$ a
continuous increasing union with $w_\delta \in \CP_\ll(\alpha)$ and $u
\sub w_0$ then there is a club $C$ so that for every $\delta \in C \cap
S'(G)$ there is $u' \in U_\delta$ so that $(w_\delta,u) \cong
(\otp(w_\delta),u')$.
\ee \ed

\noindent{\bf Remarks}:  
\be
\item Although there is the possibility of confusion between the
orders $<_\delta$ and $<$ on $\Pp_\lam$, we will use $<$ for both and
the context should usually make it clear which we mean.
\item The following will be true of the trees that we eventually
build although this property will not be needed in the proof: 
if $q \in T_\delta$ then there is a generic set $G$ with a
cofinal sequence $(g_\delta)_\delta$ so that $q =
(g_{<\alpha}^\col)\res\beta$ for some $\alpha$ and $\beta$
\item If $(g_\delta)_\delta$ and
$(g'_\delta)_\delta$ are cofinal sequences for $G$ and $G'$ then there
is a club $C$ so that if $\delta \in C$ and $\eta = \dom(g_\ld^\col)$
then $g_\ld^\col = (g'_\ld)^\col \res\eta$.  In condition
\ref{remark}, we may assume that for particular cofinal sequences, $C$
satisfies this property as well as $C \cap S(G')
\sub S(G)$.  We will often use this version of condition \ref{remark}.
\item It is important to notice the following about $p \in
T_\delta$ for which $\dom(p)$ is a limit ordinal.  If $\alpha <
\dom(p)$ then $p\res\alpha \in T_\delta$ and $p\res\alpha < p$.
Hence, any such $p$ is the limit of those elements of $T_\delta$ which
are less than it.
\ee

\bl \label{basic} ($\Dl_\lambda$)  There is a collection of basic data. \el

We leave the proof of this until the next section.  For the rest of
the paper except for the proof of Lemma \ref{basic}, we will fix a
particular choice of basic data using the notation of definition
\ref{defbasic}. 

\bd A {\em weak commitment} is a sequence $\pp = \la p^\delta : \delta <
\lam \ra$ where $p^\delta : T_\delta \rightarrow \Pp_\lam$ with the
following properties:
\be \item $p^\delta(\eta) \in \Pp_{\len(\eta)}$ (We usually write
$p^\delta_\eta$ for $p^\delta(\eta)$.)
\item if $\eta \leq \nu \in T_\delta$ then $p^\delta_\eta \geq
p^\delta_\nu\res{\len(\eta)}$. \ee

We define an order on weak commitments by $\pp \leq \qq$ if for almost
all $\delta$ (i.e. on a club), $p^\delta \leq q^\delta$ pointwise.  We say that $\qq$
is stronger than $\pp$. \ed

We will identify two weak commitments $\pp$ and $\qq$ if $\pp \leq
\qq$ and $\qq \leq \pp$.

\noindent{\bf Notation}: \quad From the fixed collection of basic data one can
extract a critical weak commitment.  Define $\up *\pp = \la \up *
p^\delta : \delta < \lam \ra$ where $\up * p^\delta_\eta = \eta$.

\bd A {\em commitment} is a weak commitment which is stronger than $\up *
\pp$. \ed

\bd Suppose $G$ is generic with a cofinal sequence $(g_\delta)_\delta$
and $\pp$ is a commitment.  We say $G$ {\em meets} $\pp$ if there is a club
$C$ in $\lambda$ so that for every $\delta \in C \cap S(G)$, $\eta_\delta  =:
(g_\ld)^\col \in T_\delta$ and there is $r_\delta \in G$ so that
$\dom(r_\delta) = \dom(g_\ld)$ and
\[ r^\col_\delta \geq p^\delta_\eta \mbox{ for all } \eta \leq \eta_\delta.\]
\ed

\noindent{\bf Remark}:
If $h: \len(\eta_\delta) \rightarrow
\dom(g_\ld)$ is an order isomorphism then in the above definition the
existence of $r_\delta$ is  equivalent to saying that
$h[p^\delta_\eta] \in G$ for all $\eta \leq \eta_\delta$.

\begin{prop} \label{basestep}

If $D_i$, $i<\lambda$, are 0-density systems and 
$g \in \Pp_\gamma$ then there is a commitment $\qq (\geq \up *\pp)$
and some $G \in
\gen(\Pp_{\gamma})$, so that:
\be \item 
$g \in G$,
\item $G$ 
meets $\qq$ and
\item if $\gamma\leq \gamma' < \lambda^+$ and
$G' \in \gen(\Pp_{\gamma'})$ 
meets $\qq$, then $G'$ meets each
$D_i$.
\ee
\end{prop}

\begin{prop} \label{successor} 

Suppose $G \in \gen(\Pp_\alpha)$ and $G$ satisfies
\be \item $ \mbox{for all } g \in \Bbb P / G, h \in \Bbb P \mbox{ there is } g'
\in \Bbb P / G \mbox{ with } g' \geq g \mbox{ and either } g' \geq h
\mbox{ or } g' \perp h$ and
\item $\Pp/G$ is $\lam^+$-uniform. \ee
For $i<\lambda$, let $D_i$ be a density system over $G$, and suppose
$g \in \Pp_\gamma/G$ where $\alpha\le\gamma < \lambda^+$ and
$\pp$ is some commitment that is met by $G$.
Then there is a commitment $\qq \geq\pp$, and some $G^* \in
\gen(\Pp_{\gamma})$, so that:
\be \item 
$G \sub G^*$, $g \in G^*$;
\item $G^*$ 
meets $\qq$;
\item \label{extension} If $\gamma\leq \gamma' < \lambda^+$ and
$G' \in \gen(\Pp_{\gamma'})$ contains $G$
and meets $\qq$, then $G'$ meets each
$D_i$.
\ee
\end{prop}

\begin{lemma}
\label{limitp} Let 
$(^\alpha\pp)_{\alpha < \kappa}$ be an
increasing sequence of commitments with
$\kappa<\lambda^+$.  Then the sequence
has a least upper bound.
\end{lemma}

\begin{notation} \quad 

With the notation of the preceding
lemma, we write
\[\bigcup_{\alpha < 
\delta}\,\up\alpha\pp \mbox{ or } {}^{\ld}\pp\]
for the least upper bound of the
commitments $\up\alpha\pp$.
\end{notation}

\begin{prop} \label{limitpG}
\quad 
Suppose $G_\alpha \in
\gen(\Pp_{\zeta_\alpha})$ meets
$^\alpha \pp$ for all $\alpha < \delta$,
where $\delta < \lambda^+$, and the
$G_\alpha$ and $^\alpha \pp$ are
increasing.  Then $G_{<\delta}$ meets
$^{<\delta} \pp$.
\end{prop}
\bigskip

By combining these results we
immediately obtain a proof of Theorem
\ref{main}.
\medskip

\noindent{\bf Proof 
of Theorem \ref{main}}:
\medskip
We define a sequence of ordinals $\zeta_\alpha$, a sequence of
commitments $\up\alpha\pp$, and a sequence of $\lambda$-generic
ideals $G_\alpha\in \gen(\Pp_{\zeta_\alpha})$, so that:
\be \item $\zeta_{<\alpha}\le\zeta_\alpha;\ \up\ld\pp\le\up\delta\pp;
\ G_\ld\includedin G_\delta \mbox{ for }
\delta<\lambda^+$
\item $G_\alpha \mbox{ 
meets the commitment }\up\alpha\pp$
\item If $\zeta_\alpha\le\beta$ 
and $G\in \gen(\Pp_\beta)$ contains
$G_\alpha$ and meets
the commitment $\up\alpha\pp$, then $G$
meets each $\alpha$-density
system $D_i$ over $G_\alpha$ proposed by Player I at
stage $\alpha$ of the genericity
game. \ee

At stage 0, Player I provides some $g_0 \in \Pp$ and at most $\lambda$
many 0-density systems.  To these Player II adds all the 0-density
systems mentioned in examples \ref{ex1} and \ref{ex2}.  We now apply
Proposition \ref{basestep}  to all these 0-density
systems and $g_0$.  
This will provide us with $G_0$, $\zeta_0$ and $\up 0\pp$.

At stage $\delta$, we will have $\zeta_\ld$, $G_\ld$,
$\up\ld\pp$ defined, and by
Proposition \ref{limitpG}\, $G_\ld$ meets
$\up\ld\pp$.  Now since $G_0 \sub G_\ld$, $G_\ld$ meets $\up 0\pp$
and hence meets each of the 0-density systems from examples \ref{ex1}
and \ref{ex2}.  It follows that $G_\ld$ satisfies the conditions on
the generic set in Proposition \ref{successor}.  By Proposition
\ref{successor}  a 
suitable choice of $\zeta_\delta$,
$G_\delta$, $\up\delta\pp$ can then be made.

Now we verify that Player II wins the
genericity game using this strategy.
By construction, $g_\alpha \in G_\alpha$ for all $\alpha$.
Suppose that $D$ is a
density
system over $G_{<\alpha}$ selected by
Player I at stage $\alpha$ of the
genericity game, and 
$\beta \geq \alpha$.  
As 
$G_\beta$ meets the commitment
$\up\beta\pp$,
$\up\beta\pp\ge\up\alpha\pp$, and
$G_\beta\contains G_\alpha$, it follows
that $G_\beta$ meets $D$.
\qed

\vfill \break

\section{Proofs}

In this section we give the proofs of
the results stated in the previous
section except the
proofs of Propositions \ref{basestep} and
\ref{successor} which are
deferred to the next section.
\bigskip

\noindent {\bf 
Proof of Lemma \ref{basic}}\medskip

The Lemma states that there is a collection of basic data.
\bigspace

Let $(v_\delta)_{\delta<\lambda}$ be an enumeration of
$\CP_{<\lambda}(\lambda)$ so that $v_\delta \sub \delta$ for all
$\delta$. Let
$V_\delta=\{v_\beta:\beta<\delta\}$ for $\delta<\lambda$.  Using
$\Dl_\lambda$ and an encoding of
\[(\lambda \times 
\lambda) \cup
(\lambda \times \CP_{<
\lambda}(\lambda))\] by $\lambda$, we
can find sets $\CR_\delta\includedin
\CP(\delta\times
\delta)$ 
and 
$\CG_\delta\includedin \CP(\delta\times
V_\delta)$, such that
$|\CR_\delta|,|\CG_\delta| < \lam$
for all $\delta<\lambda$, and for any
$R\includedin \lambda\times \lambda$,
$G\includedin \lambda\times
\CP_{<\lambda}(\lambda)$, the 
set:
\[
\{\delta:R\intersect(\delta\times \delta)\in 
\CR_\delta
\mbox{ and } G\intersect (\delta\times V_\delta)
\in \CG_\delta\}\] 
is stationary.

Before defining the basic data we
establish some notation.  For each
$\alpha<\lambda^+$, we select a
bijection $i_\alpha:\alpha\with
|\alpha|$.  For simplicity we assume
$|\alpha|=\lambda$ throughout in our
notation below.  For $\delta<\lambda$,
let $\alpha_\delta$ be the order type of
$i^{-1}_\alpha[\delta]$, and let
$\pi_{\alpha\delta}:\alpha\intersect
i^{-1}_\alpha[\delta]
\iso\alpha_\delta$, $j_{\alpha\delta}=\pi_{\alpha\delta}\circ
i^{-1}_\alpha:\delta\with\alpha_\delta$.

Let
\begin{eqnarray*}
R_\alpha & = &
\{(i_\alpha(\beta),i_\alpha(\gamma)):\beta<\gamma<\alpha\};\\
R_{\alpha\delta} & = &
R_\alpha\intersect (\delta\times
\delta)\quad (\delta<\lambda).\\
\end{eqnarray*}
Then $j_{\alpha\delta}:(\delta,
R_{\alpha\delta})\iso
(\alpha_\delta,<)$. It will be important
that $j_{\alpha\delta}$ is determined by
$R_{\alpha\delta}$.

If $G$ is a $\lambda$-generic ideal in
$\Pp_\alpha$ with $\alpha<\lambda^+$,
let:
\begin{eqnarray*}
\Gh & 
= & \{(\beta,i_\alpha[u]) : (\beta, u)
\in G \}\\
\Gh_\delta & 
= & \Gh\intersect (\delta\times
V_\delta)\\ G_\delta & = &
\{(\beta,u)\in G:
(\beta,i_\alpha[u])\in(\delta\times
V_\delta)\}\\
\Gbar_\delta &=& 
\{(\beta,\pi_{\alpha\delta}[u]): (\beta,u)\in
G_\delta\}\\
\end{eqnarray*}  

Again, we can go directly from
$\Gh_\delta$ to $\Gbar_\delta$ by
applying $j_\ad$. Observe also that
$\pi_\ad$ induces an isomorphism
$\pi^*_\ad:G_\delta\iso \Gbar_\delta$.
We are primarily interested in this
collapsing map $\pi^*_\ad$, but $\Gh$
provides a better ``encoding'' of $G$
because the sets $\Gh_\delta$ increase
with $\delta$, while the sets
$\Gbar_\delta$ do not.

Let $C(G)$ be the set of
$\delta<\lambda$ for which $G_\delta$
contains a cofinal increasing
subsequence.  Then $C(G)$ is a club in
$\lambda$. For $\delta\in C(G)$,
$\Gbar_\delta$ has a least upper bound,
which will be denoted
$\union\Gbar_\delta$.
\medskip

We are now ready to define the basic data.
For $\delta <
\lambda$ we 
define $T_\delta$ as:
\begin{eqnarray*}
&\{p\in \Pp:&
\exists 
\alpha < \lambda^+\,
\exists G 
\in \gen(\Pp_\alpha)\
\exists \gamma\!:\ 
\ \delta\in C(G),\\
&& \Gh_\delta\in \CG_\delta, R_\ad\in
\CR_\delta,
\mbox{ and 
} p = [\bigcup \Gbar_\delta] \res
\gamma. \}\\
\end{eqnarray*}
Notice that $\dom(p)$ is an ordinal for every $p \in T_\delta$.
To see that $|T_\delta|<\lambda$,
we use the fact that $\Gh_\delta, R_\ad$
together determine $\Gbar_\delta$, and
also that any $p$ in $\Pp$ has fewer
than $\lambda$ distinct restrictions.
For $p,q \in T_\delta$, define the order
by; $p \leq q$ iff $p = q \res{\dom(p)}$.

Now for $G$ a $\lambda$-generic ideal in
$\Pp_\alpha$ with $\alpha<\lambda^+$,
fix a cofinal sequence
$(g^G_i)_{i<\lambda}$ in $G$, and set:
\begin{eqnarray*}
S(G)&=&\{\delta<\lambda:[g^G_{<\delta}]^\col\in  
T_\delta\};\\ 
S'(G)&=&\{\delta < \lam: \hat G_\delta \in \CG_\delta \mbox{ and }
R_{\alpha\delta} \in \CR_\delta \};\\
U_\delta & = &\{u: \exists v\in V_\delta\,\exists R\in
\CR_\delta\,\exists \alpha<\lambda^+ \quad 
(\delta,v,R)\iso(\alpha,u,<) \}
\end{eqnarray*}
Clearly $U_\delta\includedin
\CPll(\lam)$ 
and $|U_\delta| <
\lambda$. It is also straightforward to see that $S'(G)$ is
stationary.

  Let $C_1$ be
\[\{\delta\in C(G): 
g^G_\ld=\union G_\delta 
\}\]
Then $C_1$ is a club in $\lambda$,
and if $\delta \in
S'(G) \cap C_1$ then $(g_\ld^G)^\col \in T_\delta$ so $S(G)$ is
stationary.  If $(g'_\delta)_\delta$ is any other cofinal sequence for
$G$ then there is a club 
\[C = \{ \delta : g'_\ld = g^G_\ld \} \] and for every $\delta \in C
\cap S(G)$, $(g'_\ld)^\col \in T_\delta$.

Let $G \sub G^*$ be two $\lam$-generic ideals in
$\Pp_\alpha,\Pp_{\alpha^*}$ with $S(G), S(G^*)$ determined by cofinal
sequences $(g^G_\delta)_\delta, (g^{G^*}_\delta)_\delta$ respectively.
If one considers $C = \{ \delta : g^{G^*}_\ld \res \alpha = g^G_\ld
\}$, it is easy to see that $C \cap S(G^*) \sub S(G)$.

It remains to verify the 
oracle property (\ref{universality}) of
Definition \ref{defbasic}.  We fix
$\alpha<\lambda^+$, $G$
$\lambda$-generic in $\Pp_\alpha$, $u\in
\CPll(\alpha)$, and 
we let $\alpha =
\union_{\delta< \lambda} 
w_\delta$ be a continuous increasing
union with each $|w_\delta| < \lambda$ and $u \sub w_0$.
On some club $C$, 
$\otp\,
w_\delta=\alpha_\delta$ and if $v_\alpha \sub \delta$ then $\alpha <
\delta$.  So
$(w_\delta,u)\iso(\alpha_\delta,\pi_\ad[u])$.
For $\delta\in C\intersect S'(G)$ we have
\[(\delta,i_\alpha[u],R_{\alpha\delta}) \iso
(\alpha_\delta,\pi_{\alpha\delta}[u],<) \iso (w_\delta,u,<).\]  Hence,
$\pi_{\alpha\delta}[u] \in U_\delta$. \qed

\begin{notation}
{\rm \quad In the next few results we
make systematic use of the diagonal
intersection of clubs. If
$(C_\alpha)_{\alpha<\lambda}$ 
is a sequence of clubs in $\lambda$, the
diagonal intersection is defined
correspondingly as:
\[\Diag_\alpha C_\alpha 
=  \{\delta<\lambda:
\delta\in \Intersection_{\alpha<\delta} 
C_\alpha\}.\]
The diagonal intersection of such a
sequence of clubs is again a club.}
\end{notation}

\noindent{\bf Proof 
of Lemma \ref{limitp}}\quad Let
$(^\alpha\pp)_{\alpha < \kappa}$ be an
increasing sequence of commitments with
$\kappa<\lambda^+$.  We claim that the
sequence has a least upper bound.  We
may take $\kappa$ to be a regular
cardinal, with $\kappa\le \lambda$.  We
deal with the case $\kappa=\lambda$; for
$\kappa<\lambda$ our use of a diagonal
intersection below would reduce to an
ordinary intersection.

For $\beta<\lambda$ let
$C_{\beta}$ be a club such that for all $\alpha<\beta$:
\[\up\alpha 
p^\delta\le\up\beta p^\delta \mbox{ pointwise
for } \delta\in C_{\beta} .\]

Let $C=\Diag_{\beta}
C_{\beta}$.  For $\delta\in C$ 
and  $\eta\in T_\delta$,
let
$p^\delta_\eta=\union_{\alpha<\delta}
\up\alpha p^\delta_\eta$. 
Then $\pp$ is a  commitment.  We
have $\up\alpha
\pp \le 
\pp$ since $\up \alpha p^\delta\le p^\delta$ pointwise for
$\delta\in C\sm\alpha$.

Now we will check that $\pp$ is the
least upper bound of the sequence as a
commitment 
Let $\qq$ be a
second upper bound.  Let
\[C^*_\alpha=\{\delta<\lambda: q^\delta\ge\up\alpha 
p^\delta \mbox{ pointwise}
\},\]
and let $C^*=\Diag_\alpha C^*_\alpha$.
For $\delta\in C\intersect C^*$, and
$\eta\in T_\delta$, we have:
\[q^\delta_\eta\ge\union_{\alpha<\delta}\up\alpha
p^\delta_\eta=p^\delta_\eta.\]

It follows that $\pp$ is the least
upper bound.\qed
\bigskip

We divide Proposition \ref{limitpG} into
two parts.

\begin{prop} \label{limit 
G} Suppose that $(G_i)_{i<\kappa}$ is an
increasing sequence with $G_i\in
\gen(\Pp_{\alpha_i})$, $\kappa<\lambda^+$, 
and that each $G_i$ meets a fixed
commitment $\pp$.  Then $\union_i G_i$
also meets the commitment $\pp$.
\end{prop}

\proof Let $G = \bigcup_{i<\kappa} G_i$.  By consulting the proof of
Lemma \ref{mingen}, we can see that there is a cofinal sequence
$(g_j)_{j<\lambda}$ for $G$ so that if
$g^i_j=g_j\res
\alpha_i$ then
$(g^i_j)_{j<\lambda}$ is a cofinal
sequence in $G_i$. For each $i<\lambda$
let $C_i$ be a club demonstrating that
$G_i$ meets $\pp$.  In other words, for
$\delta\in C_i\intersect S(G_i)$ we have
$r^i_\delta\in G$ with:
\be \item 
$\dom r^i_\delta=\dom g^i_\ld$, $\eta^i_\delta = (g^i_\ld)^\col \in
T_\delta$; and
\item $[r^i_\delta]^\col 
\geq p^\delta_{\eta}$ for all $\eta \leq \eta^i_\delta$.
\ee
By Definition \ref{defbasic}.5, we may also suppose that $C_i\intersect
S(G)\includedin S(G_i)$.

We consider the case  $\kappa = \lam$  (Use ordinary intersection
instead of diagonal intersection  when $\kappa < \lam$).

Let
\[C =  \Diag_i C_i \cap \{\delta<\lambda: 
\sup\dom g_\ld =\sup_{i<\delta}\alpha_i\}\]

If $\delta \in C \cap S(G)$ then we can find $r_\delta \in G$ so that
$r_\delta \geq r^i_\delta$ for $i < \delta$ and
\[ \dom(r_\delta) = \bigcup_{i<\delta} \dom(g^i_\ld) =
\bigcup_{i<\delta} \dom(g_\ld\res \alpha_i) = \dom(g_\ld).\]
Clearly, $r_\delta^\col \geq p^\delta_\eta$ for all $\eta <
(g_\ld)^\col = \eta_\delta \in T_\delta$.
However, since $p^\delta_{\eta_\delta} \res{\len(\eta)} \leq
p^\delta_\eta$ for all $\eta \leq \eta_\delta$ and
$p^\delta_{\eta_\delta} = \bigcup_{\eta < \eta_\delta}
p^\delta_{\eta_\delta}\res{\len(\eta)}$, $r_\delta^\col \geq
p^\delta_{\eta_\delta}$. \qed

\begin{prop} \label{limit 
p} Suppose $G \in \gen(\Pp_\alpha)$
meets an increasing sequence of
commitments $(\mbox{$^\gamma\pp$} : \gamma
< \kappa)$ where $\kappa <
\lambda^+$.  
Then $G$ meets $\bigcup_{\gamma <
\kappa}\up\gamma\pp$. 
\end{prop}

\proof 
\ Again, we treat only the case when $\kappa = \lam$.
Let $(g_i)_{i < \lambda}$ be a cofinal sequence in $G$.
For each $\gamma$, let $C_\gamma$ witness
the fact that $G$ meets $\up\gamma\pp$.
That is, for all $\delta \in C_\gamma
\cap S(G)$ 
there is $r_\delta^\gamma \in G$ so that:
\be \item 
$\dom r^\gamma_\delta = \dom g_\ld$, $\eta_\delta = g_\ld^\col \in
T_\delta$; and
\item $(r^\gamma_\delta)^\col \geq \up \gamma p^\delta_\eta$ for all
$\eta \leq \eta_\delta$.
\ee
We may also suppose that $C_\gamma$
witnesses the relation $\up i \pp \le
\up \gamma \pp$ for $i<\gamma$. Hence we may assume
$r^i_\delta\le r^\gamma_\delta$ for 
$\delta\in C_\gamma\intersect S(G)$.
Let $C=\Diag_\gamma C_\gamma$.

For $\delta\in C\intersect S(G)$, let $r_\delta =
\Union_{i<\delta} r^i_\delta$. 
Then on $C\intersect S(G)$, 
$\dom r_\delta =
\dom g_\ld$ and $r^\col_\delta 
= \Union_{i<\delta}[r^i_\delta]^\col
\geq\Union_{i<\delta} \up i 
p^\delta_{\eta}
= \up{<\kappa} p^\delta_{\eta}$ for all $\eta \leq \eta_\delta$.  The
last equality follows from the proof of Lemma \ref{limitp}.
\qed
  
\medskip

Proposition \ref{limitpG} is an
immediate consequence of the preceding
two propositions.

\vfill \break

\section{Proof of Proposition \protect{\ref{successor}}}

 We
recall the statement of Proposition \ref{successor}.
\bigskip

\noindent{\bf Proposition 
\ref{successor}}

\sl
Suppose $G \in \gen(\Pp_\alpha)$ and $G$ satisfies
\be \item $ \mbox{for all } g \in \Bbb P / G, h \in \Bbb P \mbox{ there is } g'
\in \Bbb P / G \mbox{ with } g' \geq g \mbox{ and either } g' \geq h
\mbox{ or } g' \perp h$ and
\item $\Pp/G$ is $\lam^+$-uniform. \ee
For $i<\lambda$, let $D_i$ be a density system over $G$, and suppose
$g \in \Pp_\gamma/G$ where $\alpha\le\gamma < \lambda^+$ and
$\pp$ is some commitment that is met by $G$.
Then there is a commitment $\qq \geq\pp$, and some $G^* \in
\gen(\Pp_{\gamma})$, so that:
\be \item 
$G \sub G^*$, $g \in G^*$;
\item $G^*$ 
meets $\qq$;
\item If $\gamma\leq \gamma' < \lambda^+$ and
$G' \in \gen(\Pp_{\gamma'})$ contains $G$
and meets $\qq$, then $G'$ meets each
$D_i$.
\ee

\rm
\medskip

We are also obliged to prove Proposition \ref{basestep} as well.  The
proof is very similar to the proof of Proposition \ref{successor} and
so we will only highlight the formal differences at the end of the proof.

\noindent{\bf Proof of Proposition \ref{successor}}:  Let $\gamma =
\bigcup_{\delta < \lam} 
w_\delta$ be a continuous increasing union with $|w_\delta| < \lam$.
Set $\gamma_\delta = \otp(w_\delta)$ and choose $\zeta_\delta$ so that
$\gamma_\delta + \zeta_\delta \geq \hgt(T_\delta)$.  Let 
\[h_\delta :
\gamma_\delta + \zeta_\delta \rightarrow w_\delta \cup [\gamma,\gamma
+ \zeta_\delta)\] be an order isomorphism.

Fix a cofinal sequence $(g_\delta)_\delta$ for $G$.  Since $G$ meets
$\pp$, there is a club $C$ so that for all $\delta \in C \cap S(G)$ we have
$g_\ld^\col =: \eta_\delta \in T_\delta$, and there is $r_\delta \in G$
so that $\dom(r_\delta) = w_\delta \cap \alpha$ and $r^\col_\delta \geq p^\delta_\eta$ for all $\eta \leq
\eta_\delta$.  We may also assume that $\dom(g_\ld) = w_\delta \cap
\alpha$ for all $\delta \in C$.  

Now we build the commitment $\qq$.  If $\delta \not \in C \cap S(G)$
then let $q^\delta = p^\delta$.  Fix then $\delta \in C \cap S(G)$.
For $i < \delta$, $\zeta \leq
\zeta_\delta$ and $u \in U_\delta$, $u \sub \gamma_\delta + \zeta_\delta$, let
\[D^\zeta_i(u) = 
\{r\in \Pp_{\gamma_\delta+\zeta}: h_\delta[r]\in
D_i(h_\delta[u],w_\delta\union
[\gamma,\gamma+\zeta)) \}\]

Let $\Pp_G[T_\delta]$ be the set of functions $\bar p: T_\delta \rightarrow
\Pp_\lam$ so that \be \item $\bar p(\eta) \in \Pp_{\len(\eta)}$ for all
$\eta \in T_\delta$,
\item if $\eta \leq \nu$ then $\bar p(\eta) \geq \bar p(\nu)
\res{\len(\eta)}$ and 
\item if $\eta$ is comparable with $\eta_\delta$ then
$h_\delta[\bar p(\eta)] \in \Pp/G$. \ee
We will write $\bar p_\eta$ for $\bar p(\eta)$.

\noindent{\bf Remark}:  Since $G$ meets $\pp$, if $\delta \in C
\cap S(G)$ then 
$p^\delta \in \Pp_G[T_\delta]$.  To see this, we must show that if
$\eta$ is comparable to $\eta_\delta$ then $h_\delta[p^\delta_\eta]
\in \Pp/G$.  If $\eta \leq \eta_\delta$ then since $\delta \in C \cap
S(G)$, $h_\delta[p^\delta_\eta] \in G$.  Suppose $\eta \geq \eta_\delta$.  Now
$p^\delta_\eta \res \len(\eta_\delta) \leq p^\delta _{\eta_\delta}$
and since $w_\delta \cap \alpha = \dom(g_\ld)$ we have
$h_\delta[p^\delta_\eta]\res \alpha \leq
h_\delta[p^\delta_{\eta_\delta}]$ so $h_\delta[p^\delta_\eta] \in \Pp/G$.

\begin{prop}\label{prop5.1}
There is a 
$q^\delta \in
\Pp_G[T_\delta]$ with $q^\delta \geq p^\delta$ pointwise and so that for
every $u \in U_\delta$, $i < \delta$ if $\eta' \in T_\delta$, $\eta'
\geq \eta_\delta$ with $\len(\eta') = \gamma_\delta + \zeta$ and $u
\sub \gamma_\delta + \zeta$ then
$q^\delta_{\eta'} \in D^\zeta_i(u)$.
\end{prop}

To obtain this $q^\delta$ we use a claim whose proof we postpone.

\begin{claim} \label{commitclaim}
 If $\bar q \in \Pp_G[T_\delta]$, $u \in U_\delta$, $i < \delta$, $\zeta
\leq \zeta_\delta$ and $\eta^* \in T_\delta$, $\eta^* \geq \eta_\delta$
with $\len(\eta^*) = \gamma_\delta + \zeta$ and $u
\sub \gamma_\delta + \zeta$  then there is $\bar r \in
\Pp_G[T_\delta]$ so that $\bar r \geq \bar q$ pointwise and $\bar r_{\eta^*}
\in D^\zeta_i(u)$. 
\end{claim}

{\bf Proof of Proposition \ref{prop5.1}}
To get the required $q^\delta$, one starts with $p^\delta$, at
limit stages take unions and at successor stages use the claim applied
to some particular $i < \delta$, $\zeta \leq \zeta_\delta$, $u \in
U_\delta$ and $\eta^* \in T_\delta$.  After at most $|U_\delta| \cdot
|\delta| \cdot |\zeta_\delta| \cdot |T_\delta|$ stages we will have
produced $q^\delta$. \qed

Now we turn to the construction of $G^*$, a $\lam$-generic ideal in
$\Pp_\gamma$ meeting $\qq$ with $G \sub G^*$ and $g \in G^*$, hence
completing the proof of Proposition \ref{successor}.

Fix an enumeration $(s_\delta)_\delta$ of $\Pp_\gamma$.  $G^*$ will be
the downward closure of an increasing sequence $(g^*_\delta)_\delta$ which is
constructed inductively starting with $g^*_0 = g$.  We shall guarantee
that $g^*_\delta \in \Pp_\gamma/G$ for each $\delta$.

At stage $\delta$, if $\delta \in C \cap S(G)$, $\dom(g^*_\ld) =
w_\delta$, $(g^*_\ld)\res\alpha = g_\ld$ and $(g^*_\ld)^\col = \eta^*
\in T_\delta$ then let $h:\len(\eta^*) \rightarrow \dom(g^*_\ld)$ be
an order isomorphism and let $\hat g_\delta \in \Pp/G$ be chosen so
that $\dom(\hat g_\delta) = \dom(g^*_\ld), \hat g_\delta \geq g^*_\ld$
and $\hat g_\delta \geq h_\delta[q^\delta_{\eta}]$ for every $\eta
\leq \eta^*$.  This can be accomplished because $\Pp/G$ is
$\lam^+$-uniform by assumption and we guaranteed that
$h_\delta[q^\delta_{\eta}] \in \Pp/G$ when we built the commitment
$\qq$.

If any of the above conditions fail, let $\hat g_\delta = g^*_\ld$.
In either case, use the assumption on $G$ to find $g^*_\delta$ so that
$g^*_\delta \in \Pp_\gamma/G$ with $\hat g_\delta \geq g^*_\ld$ and
either $g^*_\delta \geq s_\delta$ or $g^*_\delta \perp s_\delta$.

It follows easily now that $G^* \in \gen(\Pp_\gamma)$, $G \sub G^*$
and $g \in G^*$.  We now show that $G^*$ meets $\qq$.  Choose $C_1$ so
that $C_1 \cap S(G^*) \sub C \cap S(G)$ and for all $\delta \in C_1$,
$\dom(g^*_\ld) = w_\delta$ and $(g^*_\ld) \res \alpha = g_\ld$. If
$\delta \in C_1 \cap S(G^*)$ then $(g^*_\ld)^\col = \eta^* \in
T_\delta$ so from considerations at stage $\delta$, $\hat g_\delta \in
G^*$ and $(\hat g_\delta)^\col \geq q^\delta_\eta$ for all $\eta \leq
\eta^*$.  It follows that $G^*$ meets $\qq$.

Now suppose $G'$ meets $\qq$, $G' \in \gen(\Pp_{\gamma'})$ contains
$G$ with $\gamma \leq \gamma' < \lam^+$.  Fix a cofinal sequence
$(g'_\delta)_\delta$ for $G'$, a density system $D_i$ and $u \in
\CP_\ll(\gamma')$.  We want to find $w$ so that $u \sub w$ and
$D_i(u,w) \cap G' \neq \emptyset$.  Write $\gamma' \sm \gamma$ as a
continuous increasing union $\bigcup_{\delta < \lam} w'_\delta$ with
$w'_\delta \in \CP_\ll(\gamma' \sm \gamma)$.

There is a club $C_2$ with the following properties:
\be \item if $\delta \in C_2$, $(g'_\ld) \res\alpha = g_\ld$ and
$\dom(g'_\ld) = w_\delta \cup w'_\delta$;
\item $C_2 \cap S'(G') \sub C_2 \cap S(G') \sub C \cap S(G)$;
\item for $\delta \in C_2 \cap S(G')$ there is $r_\delta \in G'$ with
$\dom(g'_\ld) = \dom(r_\delta)$,
$g'_\ld \leq r_\delta$, $(g'_\ld)^\col = \eta'_\delta \in T_\delta$
and $r^\col_\delta \geq q^\delta_\eta$ for all $\eta \leq
\eta'_\delta$; and
\item for $\delta \in C_2 \cap S'(G')$ there is $u' \in U_\delta$ so
that \[ f_\delta: (\otp(w_\delta \cup w'_\delta), u') \iso (w_\delta
\cup w'_\delta,u) .\] 
\ee

This can be obtained by refering to the definition of basic data,
Definition \ref{defbasic} and Lemma \ref{basic}.  In particular,
condition 4 follows from the oracle property.

Now choose $\delta \in C_2 \cap S'(G')$ with $i < \delta$. Let
$\zeta = \otp(w'_\delta)$. 
Since $\dom(g'_\ld) = w_\delta \cup w'_\delta$ and $\eta'_\delta =
(g'_\ld)^\col \in T_\delta$, we have that $\gamma_\delta +
\zeta \leq \hgt(T_\delta)$.  Moreover, $\eta_\delta =
g^\col_\ld \in T_\delta$ and $\eta_\delta \leq \eta'_\delta$.  Hence,
\[ h_\delta[q^\delta_{\eta'_\delta}] \in D_i(h_\delta[u'],w_\delta
\cup [\gamma,\gamma + \zeta))\] and
$r^\col_\delta \geq q^\delta_{\eta'_\delta}$ with $r_\delta \in G'$.

By the indiscernibility of the density systems, we have
\[ f_\delta[q^\delta_{\eta'_\delta}] \in D_i(u,w_\delta \cup
w'_\delta)\] since 
\[ f_\delta h^{-1}_\delta : (w_\delta
\cup [\gamma,\gamma + \zeta),h_\delta[u']) \iso
(w_\delta \cup
w'_\delta,u).\]

By the indiscernibility of $\Pp$, we have $r_\delta \geq
f_\delta[q_{\eta'_\delta}^\delta]$.  Since $r_\delta \in G'$,
\[f_\delta[q^\delta_{\eta'_\delta}] \in G' \cap D_i(u,w_\delta \cup
w'_\delta)\] so $G$ meets $D_i$. \qed

It remains to prove Claim \ref{commitclaim}.

\noindent{\bf Proof of Claim \ref{commitclaim}}:  Consider the set
\[S = \{ h_\delta[\bar q_\eta] : \eta \leq \eta^* \} \]
which is a subset of $\Pp/G$.  By the compatibility
condition in the definition of $\Pp_G[T_\delta]$, $S$ is also a
compatible set so we can choose $r'_{\eta^*} \in \Pp_{\len(\eta^*)}$
so that $r'_{\eta^*} \geq \bar q_\eta$ for all $\eta \leq \eta^*$ and
$h_\delta[r'_{\eta^*}] \in \Pp/G$ since $\Pp/G$ is $\lam^+$-uniform.

Now choose $r_{\eta^*} \in D^\zeta_i(u)$ so that $r_{\eta^*} \geq
r'_{\eta^*}$. This is possible since $D_i$
is a density system over $G$.  Define
\[\bar r_\eta = 
\left\{\begin{array}{ll}
        r_{\eta^*}\res \len(\eta) & \mbox{ if
} \eta \leq \eta^* \\ \bar q_\eta &
\mbox{ otherwise. 
}
\end{array} 
\right.\] 
It is easy to check that $\bar r \in \Pp_G[T_\delta]$. \qed

To obtain a proof of \ref{basestep}, make the following changes in the
above proof.  In the statement of \ref{basestep}, there is no $G$ or
$\pp$ so at the start of the proof, one must consider all $\delta <
\lam$.  The definition of $D_i^\zeta(u)$ is the same.  We replace
$\Pp_G[T_\delta]$ with $\Pp[T_\delta]$ which is the same as
$\Pp_G[T_\delta]$ but there is no third condition.  With few formal
changes, Claim \ref{commitclaim} can be proved which allows one to
build the required $q^\delta \geq \up * p^\delta$.

The rest of the proof is almost identical except that instead of
referring to the two conditions on $G$ in the statement of Proposition
\ref{successor}, one uses the fact that $\Pp$ already possesses these
qualities by virtue of being $\lam^+$-uniform.


\newpage

\begin{thebibliography}{1}

\bibitem{Hodges}
W.~Hodges.
\newblock Building models by games.
\newblock LMSST 2.
\newblock Cambridge University Press, 1985.
     
\bibitem{MagMal}
M.~Magidor and J.Malitz.
\newblock Compact extensions of {L}({Q}).
\newblock {\em Annals of Math. Logic}, 1977.
     
\bibitem{secondI}
S.~Shelah [Sh 72].
\newblock Models with second order properties {I}: Boolean algebras with no
  undefinable automorphisms.
\newblock {\em Annals of Math. Logic \bf 14} (1978),
\newblock 57-72.
     
\bibitem{secondII}
S.~Shelah [Sh 73].
\newblock Models with second order properties {II}: On trees with no
  undefinable branches.
\newblock {\em Annals of Math. Logic} {\bf 14} (1978),
\newblock 73-87.
     
     
\bibitem{secondIII}
S.~Shelah [Sh 82].
\newblock Models with second order properties {III}: Omitting types in
  $\lambda^+$ for {L}({Q}).
\newblock {\em Archive fur Math. Logik} {\bf 21} (1981),
\newblock 1-11.
     
\bibitem{ShRub}
M.~Rubin and S.~Shelah [Sh 84].
\newblock On the elementary equivalence of automorphism groups of {B}oolean
  algebras, downward {S}kolem-{L}owenheim theorems and compactness of related
  quantifiers.
\newblock {\em Journal of Symbolic Logic \bf 45}, 1980,
\newblock 265-283.
     
\bibitem{sh7}
S.~Shelah [Sh 89].
\newblock Boolean algebras with few endomorphisms.
\newblock {\em Proc. Am. Math.} {\bf 14} (1979),
\newblock 135-142.

\bibitem{secondIV}
S.~Shelah [Sh 107].
\newblock Models with second order properties {IV}: A general method and
  eliminating diamond.
\newblock {\em Annals of Pure and Applied Logic} {\bf 25} (1983),
\newblock 183-212.
     
\bibitem{sh5}
S.~Shelah [Sh 128].
\newblock Uncountable constructions.
\newblock {\em Israel Journal of Mathematics} {\bf 51} (1985),
\newblock 273-297.
     
\bibitem{Solinac}
S.~Shelah [Sh 176].
\newblock Can you take {S}olovay's inaccessible away?
\newblock {\em Israel Journal of Mathematics} {\bf 48} (1984),
\newblock 1-47.
     
\bibitem{sh 326}
S.~Shelah [Sh 326].
\newblock{Vive la diff\'erence!}
\newblock{\ To appear.}

\end{thebibliography}
\end{document}